\documentclass[11pt]{article}
\usepackage{graphicx}
\usepackage{amsmath}
\usepackage{amssymb}
\usepackage{amsthm}
\usepackage{epsfig}
\usepackage{mathrsfs}
\usepackage{multirow}
\usepackage[usenames,dvipsnames]{color}
\usepackage{setspace}
\usepackage{multicol}
\usepackage{subcaption}
\usepackage[margin=2.54cm]{geometry}
\usepackage{lipsum}
\usepackage{bigints}
\usepackage{framed}
\usepackage{colortbl}
\usepackage{calligra}
\usepackage{mathrsfs}
\usepackage{tabularx,hhline}

\usepackage{algpseudocode,algorithm,algorithmicx}

\newtheorem{theorem}{Theorem}

\newtheorem{corollary}{Corollary}
\newtheorem{lemma}{Lemma}
\newtheorem{proposition}{Proposition}
\newtheorem{assumption}{Assumption}
\newtheorem{definition}{Definition}

\numberwithin{remark}{section}

\usepackage{subcaption}
\DeclareMathOperator{\argmin}{argmin}

\DeclareMathOperator{\st}{s.t.}

\newcommand{\X}{{\mathcal{X}}}

\def\RR{{\mathbb{R}}}
\def\PP{{\mathbb{P}}}
\def\DD{{\mathcal{D}}}
\def\Ep{{\mathbb{E}}}

\newcommand{\D}{{\mathcal{D}}}

\newcommand{\N}{{\Bbb{N}}}

\newcommand{\I}{{\Bbb I}}
\newcommand{\LL}{{\mathcal{L}}}

\usepackage{epstopdf}

\definecolor{DSgray}{cmyk}{0,1,0,0}

\title{ Distributionally Robust Optimization with 
		Confidence \\ Bands for Probability Density Functions }

\author{
Xi Chen\thanks{Stein School of Business, New York University, New York City, NY, 10012, USA. Email: {\tt xichen@nyu.edu}}
\ \ \ \ \ \ \
Qihang Lin\thanks{Tippie College of Business, University of Iowa, Iowa City, IA 52245, USA. Email: {\tt qihang-lin@uiowa.edu}}
\ \ \ \ \ \ \
Guanglin Xu\thanks{Institute for Mathematics and its Applications, University of Minnesota,
Minneapolis, MN, 55455, USA. Email: {\tt gxu@umn.edu}.}%
}

\date{ }

\begin{document}

\maketitle

\begin{abstract}

\noindent Distributionally robust optimization (DRO) has been introduced for solving 
		stochastic programs where the distribution of the random parameters is unknown 
		and must be estimated by samples from that distribution. A key 
		element of DRO is the construction of the ambiguity set, which is a set of distributions that covers the true 
		distribution with a high probability. 
		Assuming that the true distribution has a probability density function, we propose a class of ambiguity sets based on confidence bands of the true density function. The use of the confidence band enables us to take the prior knowledge of the shape of the underlying density function into consideration (e.g., unimodality or monotonicity).
		Using the confidence band constructed by density estimation techniques as the ambiguity set, we establish the convergence of the optimal 
		value of DRO to that of the stochastic program as the sample size increases. 
		However, the resulting DRO problem is computationally intractable, as it involves functional decision variables 
		as well as infinitely many constraints. To address this challenge, using the duality theory, we reformulate it into a 
		finite-dimensional stochastic program, which is amenable to a stochastic subgradient scheme as
		a solution method.  
		We compare our approach with existing state-of-the-art DRO methods on the newsvendor problem and the portfolio management problem, and the numerical results showcase the advantage of our approach. 

\mbox{}


\end{abstract}

\begin{onehalfspace}

\section{Introduction}	
	The goal of {\em stochastic programming} (SP) is to minimize the expectation 
	of an objective function that depends on both decision variables and some random parameters. 
	Assuming that the random parameters follow a distribution denoted by $P^\star$,  
	a stochastic program can be formulated as follows:
	\begin{eqnarray}
	\label{eq:spP}
	v^\star := \inf_{x \in \X}\left\{\mathbb{E}_{P^\star}[f(x,\xi)]:=\int_{\mathbb{R}^m} f(x,\xi)P^\star(d\xi)\right\},
	\end{eqnarray}
	where $x \in \RR^n$ is the vector of decision variables, $\X \subseteq \mathbb{R}^n$ is
	the feasible set,  $\xi$ is the vector of random variables taking values in $\mathbb{R}^m$ 
	with the distribution $P^\star$, and $f(x,\xi): \mathbb{R}^n \times \mathbb{R}^m \rightarrow \mathbb{R}$ is the objective function. 
	SP has been actively studied for several decades: see~\cite{birge2011introduction}, 
	\cite{shapiro2014lectures} and references thereafter.
	Suppose $\xi$ is \emph{absolutely continuous} and thus has a probability density function $p^\star$with respect to the 
	Lebesgue measure. In other words, we have $p^\star:\mathbb{R}^m\rightarrow[0,+\infty)$ such that $P^\star(A)=\int_{A}p^\star(\xi)d\xi$ for any Borel set $A\subset \mathbb{R}^m$. We can further write \eqref{eq:spP} as follows:
	\begin{eqnarray}
	\label{eq:sp}
	v^\star = \inf_{x \in \X}\left\{\mathbb{E}_{p^\star}[f(x,\xi)]:=\int_{\mathbb{R}^m} f(x,\xi)p^\star(\xi)d\xi\right\}.
	\end{eqnarray}

Despite of their popularity,~\eqref{eq:sp} is often challenging to solve since
	the distribution $P^\star$ or the density $p^\star$ of $\xi$ is rarely known in real-life applications. 
	When a set of historical observations of $\xi$ is collected, one may solve the approximation of~\eqref{eq:spP} 
	by replacing $P^\star$ with an estimated distribution from the dataset, e.g., the empirical distribution. 
	However, due to the approximation error, the decision obtained from the approximate problem may 
	be of inferior quality and thus may have an undesirable out-of-sample performance; see~\cite{bertsimas2014robust, esfahani2015data, Mak.Morton.Wood.1999}.
	An alternative approach for solving \eqref{eq:sp} with an unknown distribution is {\em distributionally robust optimization}
	(DRO), in which one constructs an {\em ambiguity set} consisting of all distributions that are likely to be $P^\star$ 
	and then minimizes the expectation of the objective function over the worst-case distribution from the ambiguity set. 
	In particular, letting $\DD$ be the ambiguity set, we can formulate DRO as follows:
	\begin{eqnarray}
	\label{eq:droP}
	\inf_{x \in \X}\sup_{P \in \DD }\left\{\mathbb{E}_{P} [f(x,\xi)] := \int_{\mathbb{R}^m} f(x,\xi) P(d\xi)\right\}.
	\end{eqnarray}
	Scarf first proposed this model for a newsvendor problem \cite{scarf1958min}, and it has
	been extensively studied over the past in operations research and operations management community since then. 

	The ambiguity sets constructed by most of the aforementioned approaches contain distributions that are not absolutely continuous.
	In fact, as shown in many of the works, the worst-case distribution in their ambiguity set corresponding to the optimal 
         decision of \eqref{eq:droP} is discrete. However, in some applications, the vector of random parameter $\xi$ is known to be 
         absolutely continuous (e.g., when $\xi$ models the price of electricity or the return of securities).	
         In this case, by solving \eqref{eq:droP} with these ambiguity sets, one may obtain a 
	solution that is hedging against a discrete distribution that will never be the true distribution. 
	This phenomenon potentially leads to over-conservative decisions in the DRO problems. 	
	
    In this paper, we consider the situation where $P^\star$ is absolutely continuous 
	and propose a family of ambiguity sets $\DD$ consisting of only absolutely continuous distributions, 
	or equivalently, distributions with density functions. With such an ambiguity set, the  DRO model corresponding to 
	\eqref{eq:sp} is given as:
	\begin{eqnarray}
	\label{eq:dro}
	\inf_{x \in \X}\sup_{p \in \DD }\left\{\mathbb{E}_{p} [f(x,\xi)] := \int_{\mathbb{R}^m} f(x,\xi) p(\xi)d\xi\right\}.
	\end{eqnarray}
Thus, the worst-case distribution corresponding to the optimal solution of \eqref{eq:dro} will be absolutely continuous.


	We note that the ambiguity set of density functions has also been considered 
	by~\cite{mevissen2013data} and~\cite{de2018distributionally}. While their ambiguity sets contain only 
	polynomial density functions, our approach does not require for similar restrictions. We also note that
	the ambiguity set considered in \cite{mevissen2013data} utilizes the kernel density estimation, as does one of our ambiguity sets. While their method must specify the Legendre polynomial series density estimator,
	our method allows for using a broader family of kernel density estimations. 
	Moreover, our ambiguity set is constructed with data samples and fully utilizes \emph{shape information} of the density function $p^\star$ (e.g., unimodality or monotonicity),
	while the ambiguity set constructed in \cite{de2018distributionally} is not based on a data-driven approach.
	Indeed, their ambiguity set is constructed using other prior knowledge on $p^\star$ (e.g., moment information).

	Other works closely related to our method include~\cite{Lam:15, LiJiang:16}.  
	Li and Jiang consider an ambiguity set with moment and
	generalized unimodal constraints~\cite{LiJiang:16}. Lam imposes
	convexity constraints on the tail of the density function~\cite{Lam:15}. 
	The main difference between our approach and theirs is that our methods use 
	the shape information on $p^\star$ and a dataset generated by $p^\star$ to construct the 
	ambiguity set, while theirs does not use data samples and directly impose shape 
	information as constraints in their optimization problems.%

	In the rest of the paper, we first propose the generic DRO in Section~\ref{sec:dro}, followed by
	the construction of our data-driven ambiguity sets by using density estimation techniques from the statistics
	literature. In particular, we will present two classes of ambiguity sets and showcase their convergence to the true 
	density function and further prove the convergence of the optimal value of~\eqref{eq:dro} to the optimal objective value of the SP in \eqref{eq:sp} as the sample 
	size increases to infinity; see detail in Section~\ref{sec:ambiguityset}. 
	The setting of our ambiguity set gives rise to a challenging problem to solve, as the resulting 
	optimization problem~\eqref{eq:dro} involves functional decision variables (i.e., the density $p$) and infinitely many constraints. 
	In Section~\ref{sec:method}, taking the exploitation of the special structure of our ambiguity set, we show that~\eqref{eq:dro} can 
	be reformulated into a finite-dimensional convex stochastic program, which is amenable to an efficient stochastic first-order method
	approach as the solution method. Finally, we validate our approach with a newsvendor problem
	and a portfolio management problem in Section~\ref{sec:experiment}. The numerical results demonstrate that our approach
	can generate decisions with superior out-of-sample performances, especially when the number of observations is limited.

	\subsection{Literature review}
	
	In the existing literature, different approaches have been utilized to construct ambiguity sets. We briefly review some popular approaches as follows:
	\begin{itemize}
		\item A moment-based ambiguity set is often constructed to consist of all distributions 
		that share common marginal or cross moments; see~\cite{bertsimas2005optimal, calafiore2006distributionally, chen2016distributionally, delage2010distributionally, de2018distributionally, dupavcova1987minimax, erdougan2006ambiguous, ghaoui2003worst, Hanasuanto:15, LiJiang:16, Natarajan:17, 
			vandenberghe2007generalized, zymler2013distributionally, zymler2013worst}. 
		DRO problems with moment-based ambiguity sets in the format of~\eqref{eq:dro} 
		can usually be reformulated into tractable conic programs, e.g., second-order cone programming problems
		or semidefinite programming problems. However, the
		constructed ambiguity sets $\DD$ are not guaranteed to converge to the true distribution 
		$P^\star$, as the size of the historical observations increases to infinity, 
		although the estimations of the moments of the random variables are 
		guaranteed to converge to their true values.
		
		\item A distance-based ambiguity set is constructed by using some distance function to measure
		the distance between two distributions in the probability space. In fact, such an ambiguity 
		set can be considered as a ball centered at a reference distribution, e.g., the empirical distribution,
		in the space of probability distributions. The distance functions considered in the literature 
		include Kullback-Leibler divergence~\cite{hu2013kullback, jiang2015data},
		$\phi$-divergence\cite{ben2013robust, Duchi:16, klabjan2013robust}, Prohorov
		metric~\cite{erdougan2006ambiguous}, empirical Burg-entropy divergence balls
		\cite{Lam:15}, and Wasserstein metric~\cite{pflug2007ambiguity, wiesemann2014distributionally, esfahani2015data, gao2016distributionally, Chen:17Wasserstein}.
		Many distance-based ambiguity sets have both asymptotic and finite-sample
		convergences. However, there is evidence that the resulting DRO problems have the
		tendency to be more challenging to solve compared to their counterparts with 
		moment-based ambiguity sets; see~\cite{esfahani2015data}.  
		
		\item More recently, hypothesis-test-based ambiguity sets have been proposed; 
		see~\cite{bertsimas2014robust, bertsimas2013data}.
		Based on a hypothesis test, e.g., $\chi^2$-test, $G$-test, etc., and a confidence level,
		these approaches construct ambiguity sets consisting of the distributions that pass the
		hypothesis test with a given set of historical data. The methods that we use in this paper 
		belong to this category.
		
		\item Likelihood approaches are also considered to construct ambiguity sets 
		in the literature; see, e.g., ~\cite{Duchi:16, LamZhou:16, qian2015composite, wang2016likelihood}.
		The likelihood approaches construct ambiguity sets consisting of all distributions
		that make a set of observations to achieve a certain level of likelihood.
	\end{itemize}

	\subsection{Notation and terminology}
	Let $\text{Proj}_{\mathcal{V}}(\cdot)$ denote the Euclidean projection operator on
	to the set $\mathcal{V}$, i.e.,
	$
	\text{Proj}_{\mathcal{V}}(u)=\argmin_{v\in \mathcal{V}}\|v-u\|_2^2.
	$
	Let $\I_{E}(\xi)$ be the indicator function that equals one when $\xi\in E$ and zero
	when $\xi\notin E$. Unless stated otherwise, the terms ``almost every'', ``measurable'',
	and ``integrable'' are defined with respect to the Lebesgue sense.
	For an extended-real valued function $g$ on $\mathbb{R}^n$, let $\text{epi}(g)$ be its
	epigraph, $\text{dom}(g)$ be its domain, $\partial g$ be its subdifferential, and
	$g'\in\partial g$ be its any subgradient. We call $g$ a multifunction if it maps a point in
	$\mathbb{R}^m$ to a subset of $\mathbb{R}^n$. A multifunction $g$ is closed valued if
	$g(\xi)$ is a closed subset of $\mathbb{R}^n$ for every $\xi\in\mathbb{R}^m$ and is
	measurable if, for every closed set $\Xi\subset\mathbb{R}^n$, the set
	$g^{-1}(\Xi):=\{\xi\in\mathbb{R}^n|g(\xi)\cap\Xi\neq\emptyset\}$ is measurable. We define
	$\text{dom}(g):=g^{-1}(\mathbb{R}^n)$. A mapping $G:\text{dom}(g)\rightarrow\mathbb{R}^n$
	is called a measurable (integrable) selection of $g$ if it is measurable (integrable) and
	$G(\xi)\in g(\xi)$ for every $\xi\in\text{dom}(g)$. We use $\int\partial g(\xi)d\xi$ to represent
	the set $\{\int\partial G(\xi)d\xi|G\text{ is an integrable selection of }g\}$.

	\section{Data-Driven Distributionally Robust Optimization} \label{sec:dro}
	In this paper, we consider an ambiguity 
	set that consists of the density functions whose value is between two known non-negative functions. To construct such a set, 
	we assume that there exists a set of $N$ 
	independent realizations of the random variable $\xi$ (i.e., samples from $p^*$) denoted by
	\[
	\widehat \Xi_N := \left\{ \hat \xi^1, \ldots, \hat \xi^N \right\} \subseteq \Xi. 
	\]
	Then, for a given $\alpha \in (0,1)$, we construct two functions $l_{\alpha}:\mathbb{R}^m\rightarrow[0,+\infty]$ and $u_{\alpha}:\mathbb{R}^m\rightarrow[0,+\infty]$ based on
	$\widehat \Xi_N$, $\alpha$, and some prior 
	knowledge on $p^*$ (e.g., its shape property) such that 
	\begin{equation}\label{eq:cov}
	\PP \left \{ l_{\alpha}(\xi) \leq p^*(\xi) \leq  u_{\alpha}(\xi), \ \forall \, \xi \in [a,b] \right \} \geq 1-\alpha.
	\end{equation}
	We call the pair $(l_{\alpha},u_{\alpha})$
	the \emph{confidence bands} for the density functions $p^*$ at a confidence level of $1-\alpha$ and $\alpha$ is called the significance level. 
	We will introduce two methods to construct such a
	ambiguity set in Section \ref{sec:ambiguityset}. 
	
	Using $(l_{\alpha},u_{\alpha})$, we can construct an ambiguity set that contains $p^*$ with a confidence level of $1- \alpha$. More specifically, let $\LL$ be the space of  all non-negative Lebesgue-measurable functions on $\RR^m$. We consider the following ambiguity set:
	\begin{eqnarray}
	\label{eq:Pset}
	\DD(\widehat \Xi_N, \, \alpha) :=\left\{ p \in \LL \big| \begin{array}{c}
	l_{\alpha}(\xi) \le p(\xi) \le u_{\alpha}(\xi), \ \forall \, \xi \in  \Xi,  
	\int_{\Xi} p(\xi) \, d\xi = 1
	\end{array} \right\},
	\end{eqnarray}
	which satisfies
	$
	\PP \left( p^*\in \D(\widehat \Xi_N, \, \alpha)\right) \geq 1 - \alpha
	$ 
	according to \eqref{eq:cov}.

With $\D(\widehat \Xi_N, \, \alpha)$, 
	we can specify the distributionally robust optimization problem 
	in~\eqref{eq:dro} as
	\begin{eqnarray} 
	\label{equ:aro}
	v_{\D(\widehat \Xi_N, \, \alpha)}^* :=\inf_{x \in \X} \sup_{p \in \DD(\widehat \Xi_N, \, \alpha)} \int_{\Xi} f(x,\xi) p(\xi) d\xi.
	\end{eqnarray}
	Immediately, we have the following result
	for the optimal value of \eqref{equ:aro}.
	\begin{proposition}
		Suppose the minimal objective value of \eqref{equ:aro} is finite and achieved at $\widehat x_N\in\X$. Let $\widehat v_N := \Ep_{p^*}[f(\widehat x_N, \, \xi)]$.
		We have 
		$
		\PP \left( v_{\D(\widehat \Xi_N, \, \alpha)}^* \geq \widehat v_N \right) \geq 1 - \alpha.
		$
	\end{proposition}   
	
	\proof{Proof.}
	Whenever $p^* \in \D(\widehat \Xi_N, \, \alpha)$, we have
	$v_{\D(\widehat \Xi_N, \, \alpha)}^* =\sup\limits_{p \in \D(\widehat \Xi_N, \, \alpha)} \Ep_{p}[f(\widehat x_N, \, \xi)] \geq \Ep_{p^*}[f(\widehat x_N, \, \xi)]=\widehat v_N $. 
	By the construction of the ambiguity set $\D(\widehat \Xi_N, \, \alpha)$, 
	we have 
	$
	\PP\left( v_{\D(\widehat \Xi_N, \, \alpha)}^* \geq \widehat v_N \right) \geq \PP\left(p^* \in \D(\widehat \Xi_N, \, \alpha) \right) \geq 1 - \alpha.
	$
	\endproof

	\section{Data-Driven Ambiguity Sets} \label{sec:ambiguityset}
	
	In this section, we present two existing methods from the statistics literature to construct
	a confidence band for a density function based on observed data. The first method
	is applicable to only univariate distributions while the second one is
	applicable to multivariate distributions. Although we only described two confidence bands as the specific instances 
	for constructing the ambiguity set in \eqref{eq:Pset}, the optimization method we propose in Section~\ref{sec:method} can be applied to DRO with any ambiguity set in the form of \eqref{eq:Pset} with other constructions of confidence bands.
	
	\subsection{Shape-restricted confidence bands}

	In this subsection, we present the method by \cite{Hengartner:95} to construct a confidence band $(l_{\alpha},u_{\alpha})$ for $p^*$ with a confidence level of $1-\alpha$ and use it to build an ambiguity set like \eqref{eq:Pset}. Although this method only applies to a univariate density function, it is able to incorporate some shape information about $p^*$ (e.g., unimodality and monotonicity) into the construction of the ambiguity set which improves the theoretical convergence rate of the set to the true density $p^*$. 
	
	We need the following assumptions in this subsection.
	\begin{assumption}[For shape-restricted confidence bands.]
		\label{assume:Hengartner}We assume:
		\begin{itemize}			 
			\item[A1.] $\Xi$ (the support of $p^*$) is connected and contained in $[a,b]$ with known $a$ and $b$ satisfying $-\infty<a<b<+\infty$.
			\item[A2.] $p^*$ is unimodal with a known mode $\mu\in[a,b]$, meaning that $p^*$ is monotonically increasing on $[a,\mu]$ and decreasing on $[\mu,b]$. 
			\item[A3.] There exists a known constant $U$ such that $p(\mu) \leq U$ for any $\xi\in\Xi$.  
		\end{itemize}
	\end{assumption}

	In Assumption~\ref{assume:Hengartner} [A1.], we assume that $a$ and $b$ are
	finite for the simplicity of the notations in the derivation below. In fact, our results can be generalized when $a=-\infty$ and $b=+\infty$. 
	Moreover, for most applications, a conservative estimation of the range of $\xi$ is usually available, which can be directly used as $[a,b]$.    
	We also note that Assumption~\ref{assume:Hengartner} [A2.] also covers the case where $p^*$ is known to be monotonically increasing ($\mu=b$) or
	decreasing ($\mu=a$). In some real applications, the unimodality or monotonicity of a random parameter is a well-known fact and thus should be incorporated into the problem formulation. Note that the method by \cite{Hengartner:95} can also used to construct a confidence band even when the mode is only known to be in an interval $[\mu^+, \mu^-]\subset[a,b]$. 

	First, let $(\hat \xi_{(1)}, \ldots, \hat \xi_{(N)})$ be the order statistics of $\xi$ constructed from $\widehat \Xi_N$ satisfying $\hat \xi_{(1)}< \cdots<\hat \xi_{(N)}$. We choose a group size $K$ satisfying $0<K<N$ and define $M':=\lfloor N/K \rfloor$ and $M:=\lceil N/K \rceil$. Then we partition the sorted sequence $(\hat \xi_{(1)}, \ldots, \hat \xi_{(N)})$ into $M$ groups with the first $M'$ groups of size $K$ and the $M$th group (if $M'=M-1$) of size $N-KM'$. We then define $k_i:=(i-1)K+1$ for $i=1,2, \dots, M'$ and $k_i=N$ for $i=M$ if $M\neq M'$. Let $F^*:\mathbb{R}\rightarrow[0,1]$ be the cumulative density function of $\xi$. It is well-known that (citation?) the random variable $ \Delta_i := F^*(\hat \xi_{(k_i)})-F^*(\hat \xi_{(k_{(i-1)})})$ has the same distribution as
	\begin{equation}\label{eq:gamma_mc}
	\tilde\Delta_i := 
	\left\{\begin{array}{ll}\frac{\Gamma(K, 1)}{\sum_{j=1}^{M'} \Gamma(K, 1) + \Gamma(M-M'+1,1)}&\text{ if }i=1,2, \dots, M'\\\frac{\Gamma(N-KM', 1)}{\sum_{j=1}^{M'} \Gamma(K, 1) + \Gamma(M-M'+1,1)}&\text{ if }i=M\neq M'\end{array}\right.,
	\end{equation}
	where  $\Gamma(A, B)$ represents a gamma random variable with a shape parameter $A$ and a rate parameter $B$.
	Let $c^-(\alpha)$ and $c^+(\alpha)$ be two constants that satisfy
	\[
	\PP \left\{c^{-} (\alpha) \leq \Delta_i \leq c^+ (\alpha)\right\}=\PP \left\{c^{-} (\alpha) \leq \tilde\Delta_i \leq c^+ (\alpha)\right\} \geq 1-\alpha.
	\]
	Both $c^-(\alpha)$ and $c^+(\alpha)$ can be estimated to arbitrary
	precision by sampling according to \eqref{eq:gamma_mc}.  Let $\LL_\mu$ be the set of all density functions on $\RR$ with mode at $\mu$, i.e.
	$$
	\mathcal{L}_{\mu}:=\left\{p\in\LL|p(\xi)\geq 0, \int_{\RR} p(\xi)d\xi=1, \text{the mode of }p\text{ is }\mu.\right\}
	$$
	and 
	$$
	\mathcal{D}_{\mu}(\widehat \Xi,\alpha):=\left\{p\in\LL_\mu\bigg|c^{-} (\alpha) \leq\int_{\hat \xi_{(k_{i-1})}}^{\hat \xi_{(k_i)}}p(\xi)d\xi\leq c^+ (\alpha),i=2,3,\dots,M\right\}.
	$$
	Then, by the definitions of $c^-(\alpha)$ and $c^+(\alpha)$, we have 
	\[
	\PP \left\{p^*\in\mathcal{D}_{\mu}(\widehat \Xi,\alpha)\right\} \geq 1-\alpha.
	\]
	Given the property above, $\mathcal{D}_{\mu}(\widehat \Xi,\alpha)$ can be used as an uncertainty set $\mathcal{D}(\widehat \Xi,\alpha)$ in \eqref{equ:aro}. 
	However, this uncertainty set may be too large to ensure a good solution from solving \eqref{equ:aro}. Therefore, we need to further refine it using the prior knowledge about the shape of $p^*$ given in Assumption~\ref{assume:Hengartner}. Next, we show how to construct a shape-constrained confidence band using the technique from in \cite{Hengartner:95}.

	Given any $\xi \in [a,b]$, let $\widehat N$ ($M\leq \widehat N \leq M+4$) be the number of distinct elements in the set $\{a, b,\mu, \xi, \hat \xi_{(k_1)}, \ldots, \hat \xi_{(k_M)}\}$ and $z_j$ be the $j{\text{th}}$ smallest element in this set. 
	Then, we consider the following two sets of non-negative step functions on $\mathbb{R}$ defined as 
	\begin{eqnarray} \label{equ:shape_band1set}
	\mathcal{D}_{\mu}^-(\widehat\Xi,\xi)&:=&\left\{p\in\LL\left|\begin{array}{lll} p(\cdot) & = & \sum\limits_{\{j:z_{j+1}<\mu\}} \beta_j\I_{(z_j, z_{j+1}]}(\cdot)+ \max\limits_{\{j:\mu\in[z_j,z_{j+1}]\}}\beta_j \I_\mu(\cdot)  \\
	& & + \sum\limits_{\{j:\mu\in[z_j,z_{j+1}]\}}\beta_j\I_{(z_j, z_{j+1})}(\cdot) + \sum\limits_{\{j:z_j>\mu\}} \beta_{j} \I_{[z_j, z_{j+1})}(\cdot),
	\\&&\text{where }\beta_j\in[0,U]\text{ for }j=1,2,\dots,\hat N-1. 
	\end{array}\right.
	\right\}
	\end{eqnarray}
	\begin{eqnarray} \label{equ:shape_band2set}
	\mathcal{D}_{\mu}^+(\widehat\Xi,\xi)&:=&\left\{p\in\LL\Bigg|\begin{array}{lll} p(\cdot) & = &  \sum\limits_{\{j:z_{j+1}\leq\mu\}}  \beta_j\I_{[z_j, z_{j+1})}(\cdot)+ \infty\cdot \I_\mu(\xi) +\sum\limits_{\{j:z_j>\mu\}}\beta_{j} \I_{(z_j, z_{j+1}]}(\cdot),
	\\&&\text{where }\beta_j\in[0,U]\text{ for }j=1,2,\dots,\hat N-1. 
	\end{array}
	\right\}
	\end{eqnarray}
	where $\mathbb{I}(\cdot)$ denotes the indicator function. 
	Both $\mathcal{D}_{\mu}^-(\widehat\Xi,\xi)$ and $\mathcal{D}_{\mu}^+(\widehat\Xi,\xi)$ consist of piecewise constant density functions on $[a,b]$ with break points on $\{a, b,\mu, \xi, \hat \xi_{(k_1)}, \ldots, \hat \xi_{(k_M)}\}$ but with different continuous properties (i.e., left-continuity or right-continuity) on these break points. Then, the confidence band for $p^*$ given by \cite{Hengartner:95} is a pair of functions $(l_{\alpha},u_{\alpha})$ defined as
	\begin{eqnarray}\label{equ:shape_band1}
	l_{\alpha}^{SR}(\xi):=\inf\limits_{p\in\mathcal{D}_{\mu}^-(\widehat\Xi,\xi)\cap \mathcal{D}_{\mu}(\widehat \Xi,\alpha)} p(\xi)
	\quad\text{ and }\quad
	u_{\alpha}^{SR}(\xi):=\sup\limits_{p\in\mathcal{D}_{\mu}^+(\widehat\Xi,\xi)\cap \mathcal{D}_{\mu}(\widehat \Xi,\alpha)} p(\xi).
	\end{eqnarray}
	Here, the superscript ``SR'' refers to ``shape restricted''. By its definition, the value of function $l_{\alpha}^{SR}(\xi)$ ($u_{\alpha}^{SR}(\xi)$) equals the smallest (largest) value at $\xi$ of all density functions which have the form in \eqref{equ:shape_band1set} (\eqref{equ:shape_band2set}) and satisfy $c^{-} (\alpha) \leq\int_{\hat \xi_{(k_{i-1})}}^{\hat \xi_{(k_i)}}p(\xi)d\xi\leq c^+ (\alpha)$ for all $i$. 
	Then, our shape-restricted uncertainty set are defined as
	\begin{eqnarray}
	\label{eq:PsetSR}
	\DD^{SR}(\widehat \Xi_N, \, \alpha) :=\left\{ p \in \LL \big| \begin{array}{c}
	l_{\alpha}^{SR}(\xi) \le p(\xi) \le u_{\alpha}^{SR}(\xi), \ \forall \, \xi \in  \Xi,  
	\int_{\Xi} p(\xi) \, d\xi = 1
	\end{array} \right\}.
	\end{eqnarray}
	The following theorem is given by~\cite{Hengartner:95}.
	\begin{theorem}
		\label{thm:SRcover}
		$\PP \left\{p^*\in\mathcal{D}^{SR}(\widehat \Xi_N,\alpha)\right\} \geq 1-\alpha$.
	\end{theorem}
	
	We then describe the method for compute $l_{\alpha}^{SR}(\xi)$ and $u_{\alpha}^{SR}(\xi)$ for any $\xi$. 
	Given $\xi$, recall that $z_j$ denotes  the $j{\text{th}}$ smallest element in $\{a, b,\mu, \xi, \hat \xi_{(k_1)}, \ldots, \hat \xi_{(k_M)}\}$. Let $\widehat j$ and $\widetilde j$ be the indexes in $\{1,2,\dots,\widehat N\}$ such that $z_{\widehat j}=\xi$ and $z_{\widetilde j}=\mu$. Using these indexes, we define a polyhedron 
	\begin{eqnarray}\label{eq:lower_beta}
	\mathcal{H}(\widehat \Xi_N,\alpha)=\left\{\beta\in\mathbb{R}^{\widehat N}\left|
	\begin{array}{c}
	\beta_1 \leq \beta_2 \leq \cdots \leq \beta_{\widetilde j-1},\quad  \beta_{\widetilde j} \geq \beta_{\widetilde j+1} \geq \cdots \geq \beta_{\widehat N-1} \\\label{eq:lower_betac2}
	c^{-}(\alpha) \leq  \sum_{j: \xi_{(k_{i-1})} \leq z_j < \xi_{(k_{i})}}   \beta_j (z_{j+1}-z_j) \leq c^{+}(\alpha), \quad \forall \, i=2, \ldots, M   \\\label{eq:lower_betac3}
	\sum_{j=1}^{\widehat N-1} \beta_j (z_{j+1}-z_j) =1,  \\
	0\leq\beta_j\leq U,  \quad \forall \,1 \leq j \leq \widehat N-1,  
	\end{array}\right.
	\right\}.
	\end{eqnarray}
	According to \eqref{equ:shape_band1set}, \eqref{equ:shape_band2set} and \eqref{equ:shape_band1}, for any $\xi$, 
	$l_{\alpha}^{SR}(\xi)$ and $u_{\alpha}^{SR}(\xi)$
	can be computed respectively by solving the following linear programs
	\begin{eqnarray}\label{eq:beta}
	l_{\alpha}^{SR}(\xi)= \min\limits_{\beta\in \mathcal{H}(\widehat \Xi_N,\alpha)}  \beta_{\widehat j}\quad \text{ and }\quad u_{\alpha}^{SR}(\xi)= \min\limits_{\beta\in \mathcal{H}(\widehat \Xi_N,\alpha)} \beta_{\widehat j}.
	\end{eqnarray}
	Note that the first line of constraints in \eqref{eq:lower_beta} is because $\mu$ is the mode of any $p$ in $\mathcal{D}_{\mu}^-(\widehat\Xi,\xi)\cap \mathcal{D}_{\mu}(\widehat \Xi,\alpha)$ and $\mathcal{D}_{\mu}^+(\widehat\Xi,\xi)\cap \mathcal{D}_{\mu}(\widehat \Xi,\alpha)$. The second line of the constraints in \eqref{eq:lower_beta}
	corresponds to the condition $c^{-} (\alpha) \leq\int_{\hat \xi_{(k_{i-1})}}^{\hat \xi_{(k_i)}}p(\xi)d\xi\leq c^+ (\alpha)$ when $p$ is the step function in \eqref{equ:shape_band1set} and  \eqref{equ:shape_band2set}. The third line requires  that $p$ must a density function while the last line requires  $p$ is non-negative and no more than $U$ according the prior information of $p*$.
	
	We summary this procedure by \cite{Hengartner:95} for constructing a confidence band for a unimodal density function in Algorithm~\ref{alg:cbHengartner}. 
	\begin{algorithm}[h!]
		\caption{Shape-restricted confidence band $(l^{SR}_\alpha(\xi), u^{SR}_\alpha(\xi))$ at $\xi\in\Xi$}\label{alg:cbHengartner}
		\begin{algorithmic}[1]
			\Require Data $\widehat \Xi_N = \{ \hat \xi^1, \ldots, \hat \xi^N \}$ sampled from $p^*$, an interval $[a,b]$ containing $\Xi$, the mode $\mu$ of $p^*$, a constant $U\geq p^*$ on $\Xi$, a significance level $\alpha\in(0,1)$, a group size $K$ with $0<K<N$ and a point  $\xi\in\Xi$.	
			\State Let $M':=\lfloor N/K \rfloor$ and $M:=\lceil N/K \rceil$
			\State Let $\hat \xi_{(1)}, \ldots, \hat \xi_{(N)}$ be the order statistics of $\widehat \Xi_N$ with $\hat \xi_{(1)}< \cdots<\hat \xi_{(N)}$.
			\State Let $k_i:=\left\{\begin{array}{ll}(i-1)K+1&\text{ if }i=1,2, \dots, M'\\N&\text{ if }i=M\neq M'\end{array}\right.$ 
			\State Define random variables 
			\begin{equation*}
			\tilde\Delta_i :=
			\left\{\begin{array}{ll}\frac{\Gamma(K, 1)}{\sum_{j=1}^{M'} \Gamma(K, 1) + \Gamma(M-M'+1,1)}&\text{ if }i=1,2, \dots, M'\\\frac{\Gamma(N-KM', 1)}{\sum_{j=1}^{M'} \Gamma(K, 1) + \Gamma(M-M'+1,1)}&\text{ if }i=M\neq M'\end{array}\right.
			\end{equation*}
			where $\Gamma(\cdot, \cdot)$ is a gamma random variable parameterized by a shape parameter and a scale parameter.
			\State Find constants $c^-(\alpha)$ and $c^+(\alpha)$ that satisfy
			$
			\PP \left\{c^{-} (\alpha) \leq \tilde\Delta_i \leq c^+ (\alpha),i=1,\dots,M\right\} \geq 1-\alpha.
			$
			\State  Let $z_j$ be the $j{\text{th}}$ smallest value in the set $\{a, \mu, \xi, \hat \xi_{(k_1)}, \ldots, \hat \xi_{(k_M)}, b \}$ for $j=1,\dots,\widehat N$ where $\widehat N$ represents the number of distinct elements in that set.  
			\State Let $\widehat j$ and $\widetilde j$ be the indexes in $\{1,2,\dots,\widehat N\}$ such that $z_{\widehat j}=\xi$ and $z_{\widetilde j}=\mu$. Define the polyhedron $\mathcal{H}(\widehat \Xi_N,\alpha)$ in \eqref{eq:lower_beta}.
			\State Compute $(l_{\alpha}^{SR}(\xi),u_{\alpha}^{SR}(\xi))$ by solving the two linear programs in \eqref{eq:beta}.
			\Ensure $(l_{\alpha}^{SR}(\xi),u_{\alpha}^{SR}(\xi))$
		\end{algorithmic}
	\end{algorithm}

It is worthwhile to note that if the density is monotone, for example, non-decreasing, we have $\mu=b$ and the first line of the constraints 
	in \eqref{eq:lower_beta}  becomes
	$\beta_1 \leq \beta_2 \leq \cdots \leq \cdots \leq \beta_{\widehat N-1}$.
	Note also that, in the original work \cite{Hengartner:95}, no upper bound $U$ for $p^*$ is needed to construct this uncertainty set $\DD^{SR}(\widehat \Xi_N, \, \alpha) $. In this paper, we require knowing $U$ so that we can include a constraint $\beta_j \leq U$ in \eqref{eq:lower_beta}, this modification does not change the statistical property of $\DD^{SR}(\widehat \Xi_N, \, \alpha) $ (i.e., Theorem~\ref{thm:SRcover}) once we assume $p^*\leq U$. We need the constraint $\beta_j \leq U$ only to ensure $l_{\alpha}^{SR}(\xi)\leq U$ and $u_{\alpha}^{SR}(\xi)\leq U$ which are needed in our theoretical analysis.

	
	\cite{Hengartner:95} established the convergence rate
	of the constructed confidence band in the following proposition.\footnote{The original Corollary 7.2~\cite{Hengartner:95} was stated a little differently form from 	
		Proposition~\ref{lem:rate}. In fact, the conclusion there replaces $\lim$ in \eqref{eq:rateband} by $\liminf$ and equality in \eqref{eq:rateband} by $>0$. Their results were obtained by choosing the parameter $\tau=2+2\rho$ appearing in the proof of Corollary 7.2~\cite{Hengartner:95}. However, the same proof works for Proposition~\ref{lem:rate} and thus implies \eqref{eq:rateband} if we choose  $\tau=3+2\rho$. }
	
	\begin{proposition}[Corollary 7.2~\cite{Hengartner:95}] \label{lem:rate} 
		Let $\xi \in \Xi$ but $\xi \neq \mu$. Suppose $p^*$ is $(C,\rho)$-Holder continuous for constants $C>0$ and $\rho>0$, i.e., 
		$\left|p^*(\xi')- p^*(\xi)\right| \leq C \left|\xi' -\xi \right|^{\rho}$ for any $\xi$ and $\xi'$.
		Suppose $K =\left\lceil B\left(N^{2\rho}\log N\right)^{1/(1+2\rho)}\right\rceil$ for a constant $B>0$ and sufficiently large $N$. 
		We have
		\small
		\begin{equation}
		\label{eq:rateband}
		\lim_{N \rightarrow \infty} \PP \left\{ \left |u_{\alpha}^{SR}(\xi) - l_{\alpha}^{SR}(\xi) \right| \leq 4p^*(\xi) 
		\left(\sqrt{\frac{3+2\rho}{1+2\rho}}B^{-1/2}+C(p^*(\xi))^{-(\rho+1)}B^\rho\right)
		\left(\frac{\log N}{N}\right)^{\rho/(1 + 2\rho)} \right\} =1.
		\end{equation}
		\normalsize
	\end{proposition}
	
	According to the lower bound in \cite{Khasminskii:76}, the convergence rate $O\left(\left(\frac{\log N}{N}\right)^{\rho/(1 + 2\rho)}\right)$ in \eqref{eq:rateband}
	attains the minimax rate and thus is optimal (upto a constant factor) for the confidence band for a unimodal density. 
	It is also worthwhile to note that the rate
	of convergence from Kolmogorov-Smirnov distance is slower as compared to this approach. 
	As shown in \cite{Hartigan:85}, the confidence band formed
	by the Kolmogorov-Smirnov distance is only $O\left((\frac{1}{N})^{1/4}\right)$,
	which is slower than the rate of $O\left((\frac{\log N}{N})^{1/3}\right)$ in 
	Proposition \ref{lem:rate} when $\rho=1$. 
	
	Proposition~\ref{lem:rate} shows the pointwise convergence of the confidence band to the true distribution $p^*$. In the next theorem, we show that, under additional assumption, we can characterize the convergence in probability  of 
	$v_{\D^{SR}(\widehat \Xi_N, \, \alpha)}^*$ define in \eqref{equ:aro} to $v^*$ in \eqref{eq:sp}.
	
	\begin{theorem} \label{thm:convergence}
		Suppose $p^*$ is $(C,\rho)$-Holder continuous for constants $C>0$ and $\rho>0$, i.e., 
		$\left|p^*(\xi')- p^*(\xi)\right| \leq C \left|\xi' -\xi \right|^{\rho}$ for any $\xi$ and $\xi'$. Moreover, suppose $\max_{x\in\mathcal{X},\xi\in \Xi}|f(x,\xi)|<+\infty$. For any $\epsilon>0$ and $\theta\in (0,1)$, there exists an $N_{\epsilon,\theta}$ such that for any $N\geq N_{\epsilon,\theta}$, we have
		\small
		\begin{eqnarray*}
			\PP\left(\sup_{x\in\mathcal{X}}\left|\int_{\Xi} f(x,\xi)p^*(\xi)d\xi-\sup_{p \in \DD^{SR}(\widehat \Xi_N, \, \alpha)} \int_{\Xi} f(x,\xi) p(\xi) d\xi\right|\leq \epsilon\right)\geq 1-\theta
		\end{eqnarray*}
		\normalsize
		and
		$
		\PP\left(\left|v_{\D^{SR}(\widehat \Xi_N, \, \alpha)}^*-v^*\right|\leq \epsilon\right)\geq 1-\theta.
		$
	\end{theorem}
	\proof{Proof.}
	See Appendix~\ref{sec:proof2}.
	\endproof

\begin{figure}
		\caption{The data-generating distribution is $250\times\text{Beta}(5,2)$. The first row contains three subfigures with the same sample size
			($N=100$) but different significance levels $\alpha$. The second row contains three subfigures with the same significance level ($\alpha = 0.2$) but different
			sample sizes.}
		\begin{subfigure}[b]{0.33\textwidth}
			\includegraphics[width=\linewidth]{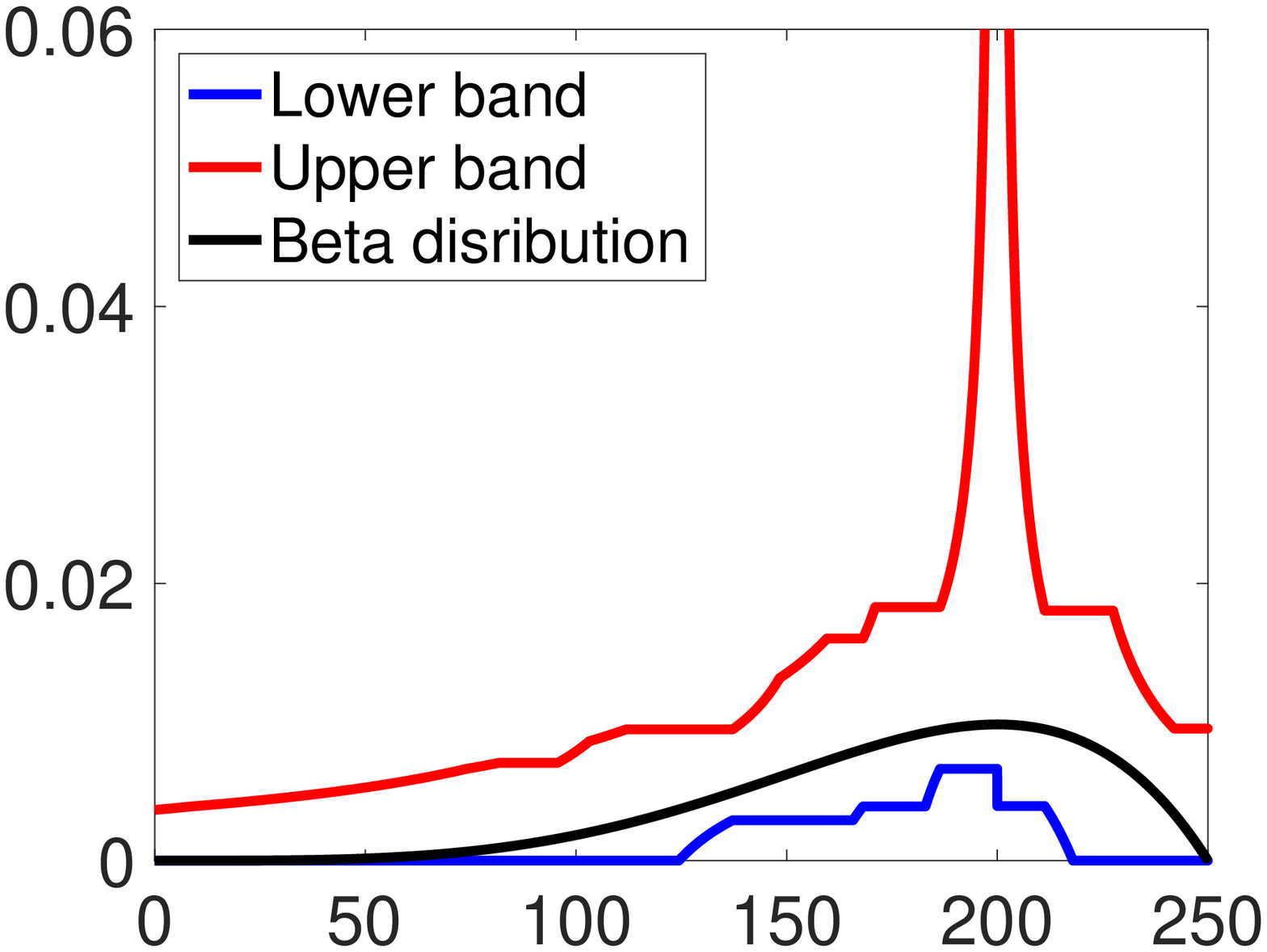}
			\caption{$\alpha=0.1$} \label{Fig:m}
		\end{subfigure}%
		\begin{subfigure}[b]{0.33\textwidth}
			\includegraphics[width=\linewidth]{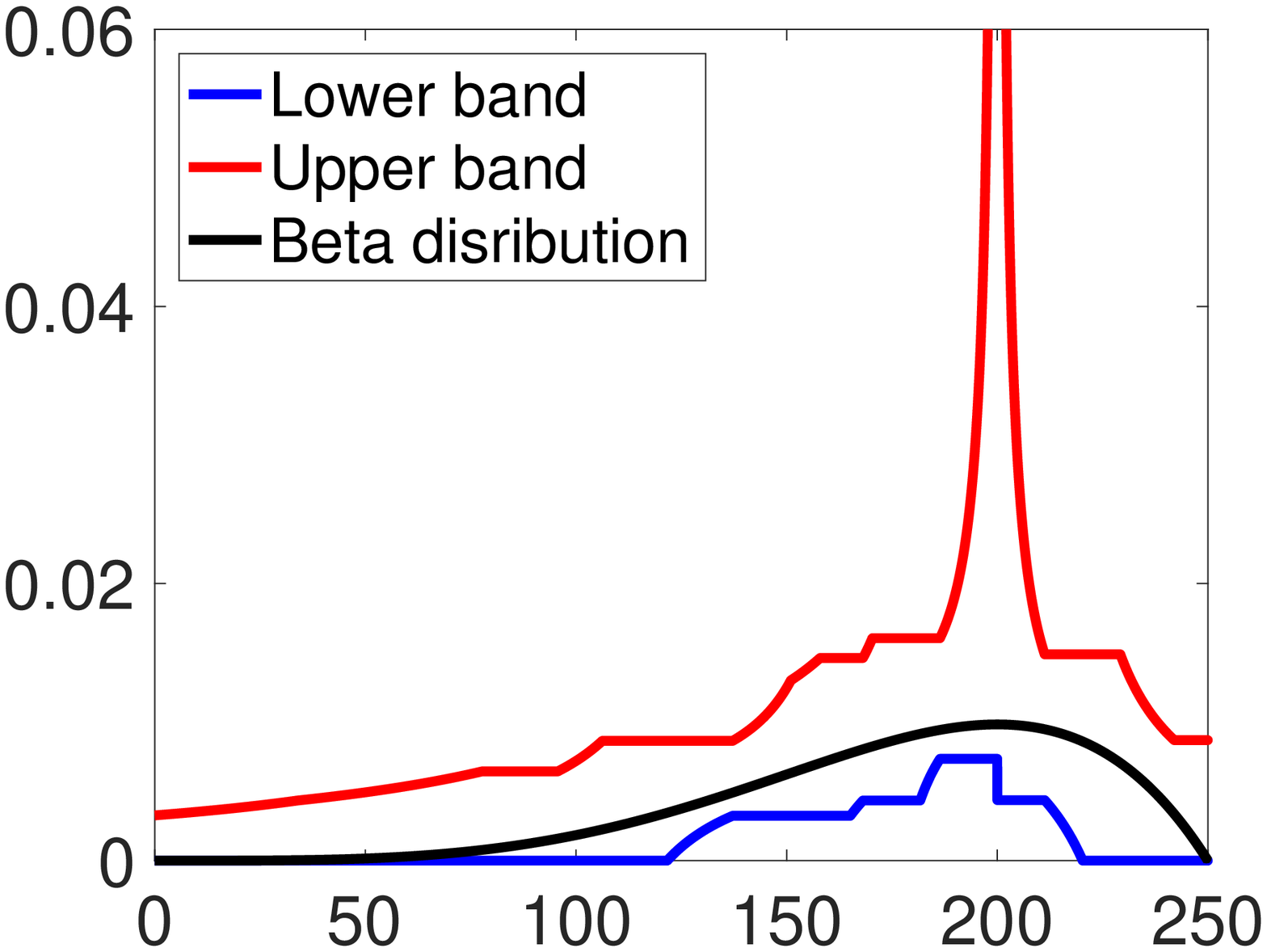}
			\caption{$\alpha=0.2$}
			\label{Fig:2m}
		\end{subfigure}%
		\begin{subfigure}[b]{0.33\textwidth}
			\includegraphics[width=\linewidth]{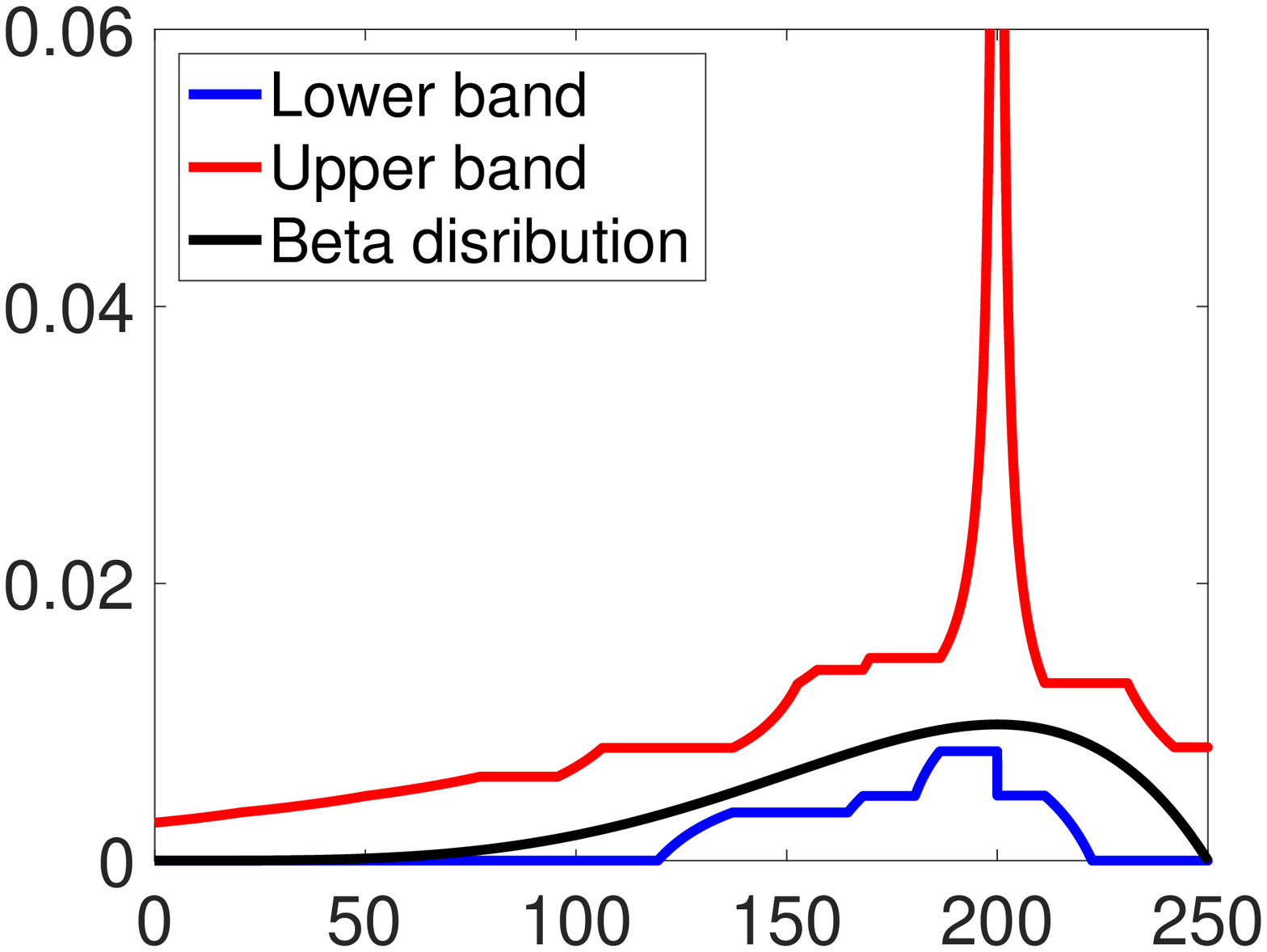}
			\caption{$\alpha=0.3$}
			\label{Fig:5m}
		\end{subfigure}%
		\quad
		\begin{subfigure}[b]{0.33\textwidth}
			\includegraphics[width=\linewidth]{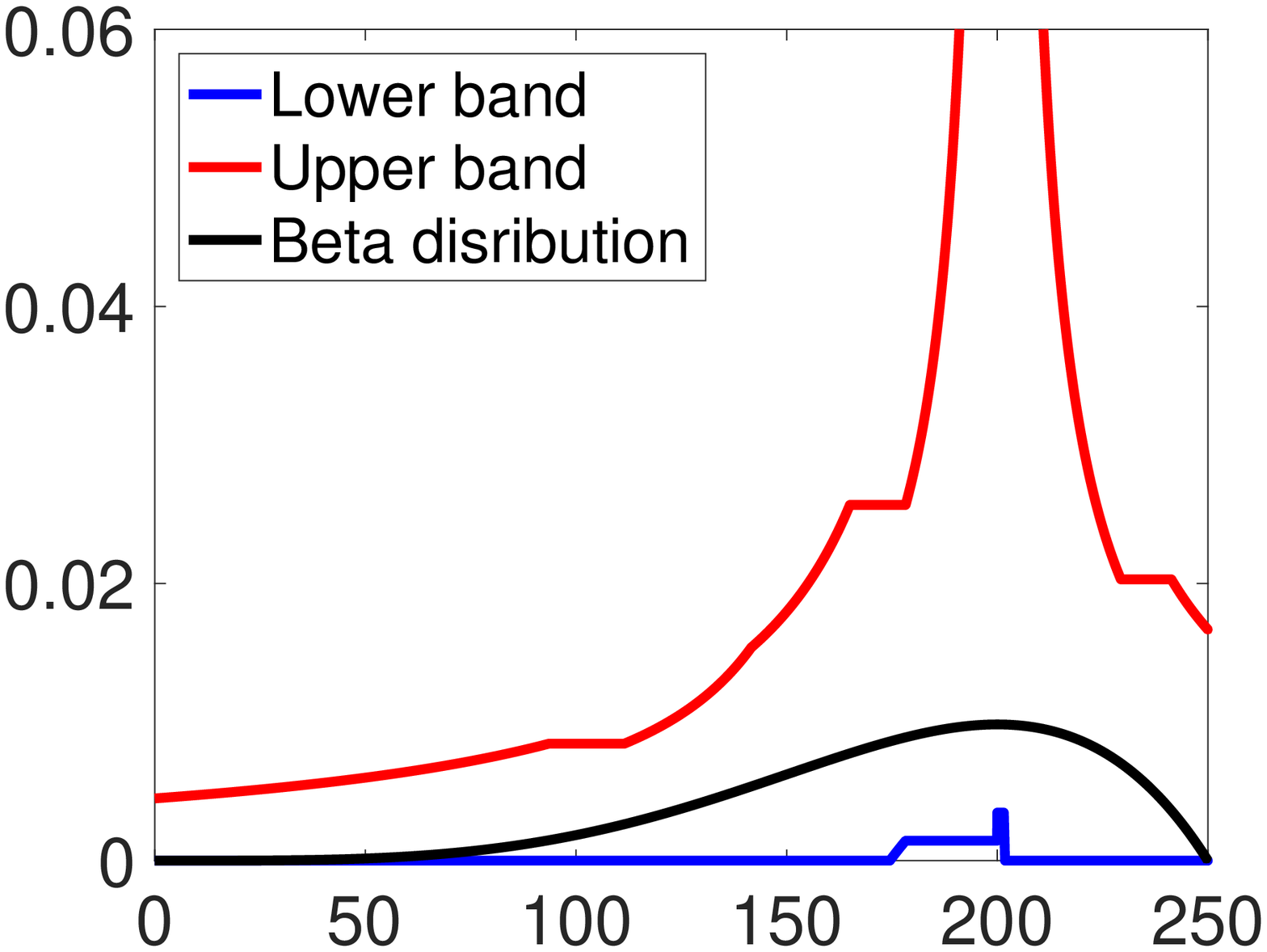}
			\caption{$N=10$}
			\label{Fig:m}
		\end{subfigure}%
		\begin{subfigure}[b]{0.33\textwidth}
			\includegraphics[width=\linewidth]{beta_100_02.eps}
			\caption{$N=100$}
			\label{Fig:2m}
		\end{subfigure}%
		\begin{subfigure}[b]{0.33\textwidth}
			\includegraphics[width=\linewidth]{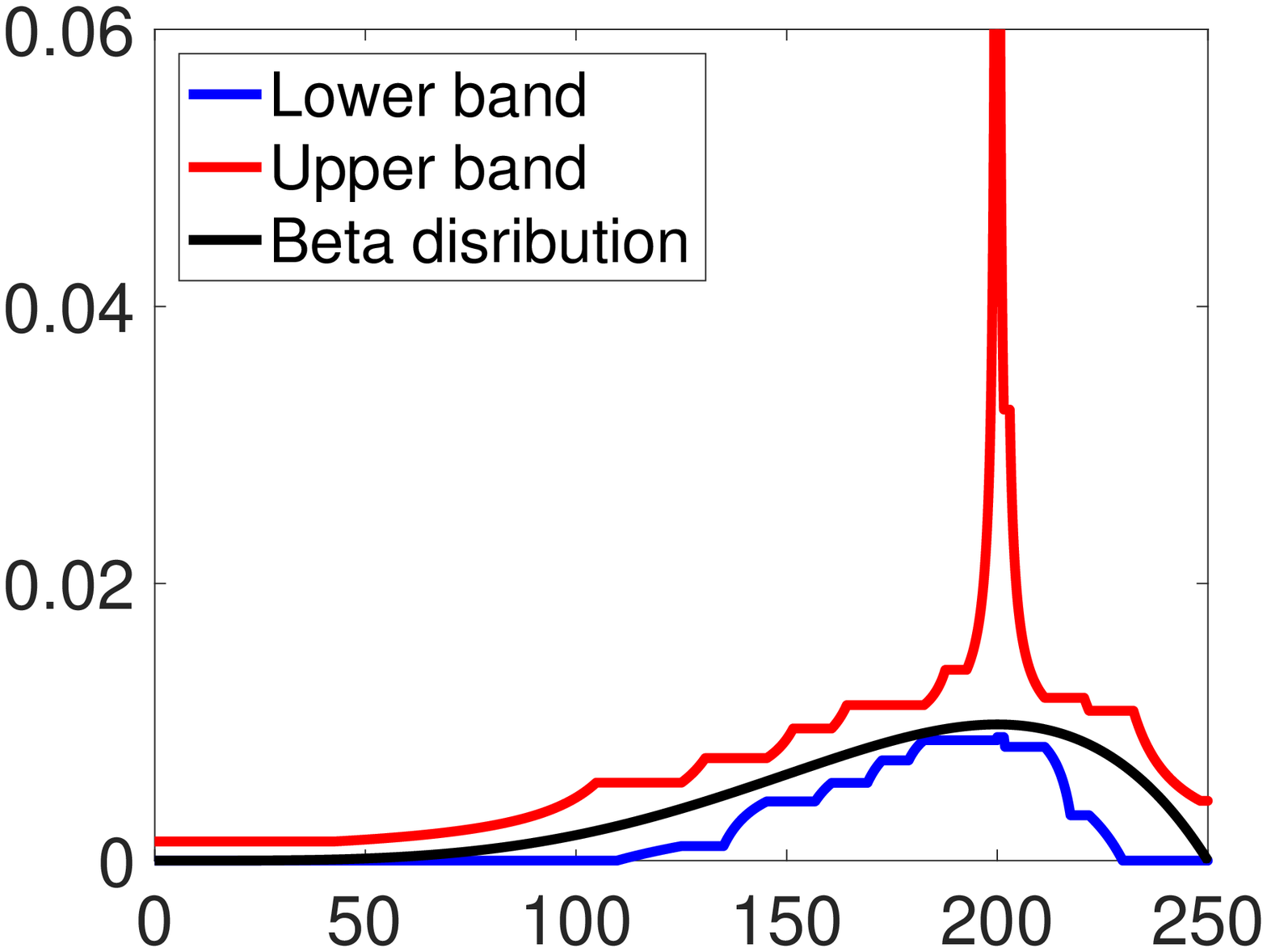}
			\caption{$N=1000$}
			\label{Fig:5m}
		\end{subfigure}
		\label{Fig:beta_mode}
	\end{figure}

\begin{figure}
		\caption{ The data-generating distribution is a trucated exponential distribution, $\text{Exp}(1/100)$
		with the support $[0,250]$. 
		          The first row contains three subfigures with the same sample size
			($N=100$) but different significance levels $\alpha$. The
			second row contains three subfigures with the same significance level ($\alpha = 0.2$)
			but different sample sizes.}
		\begin{subfigure}[b]{0.33\textwidth}
			\includegraphics[width=\linewidth]{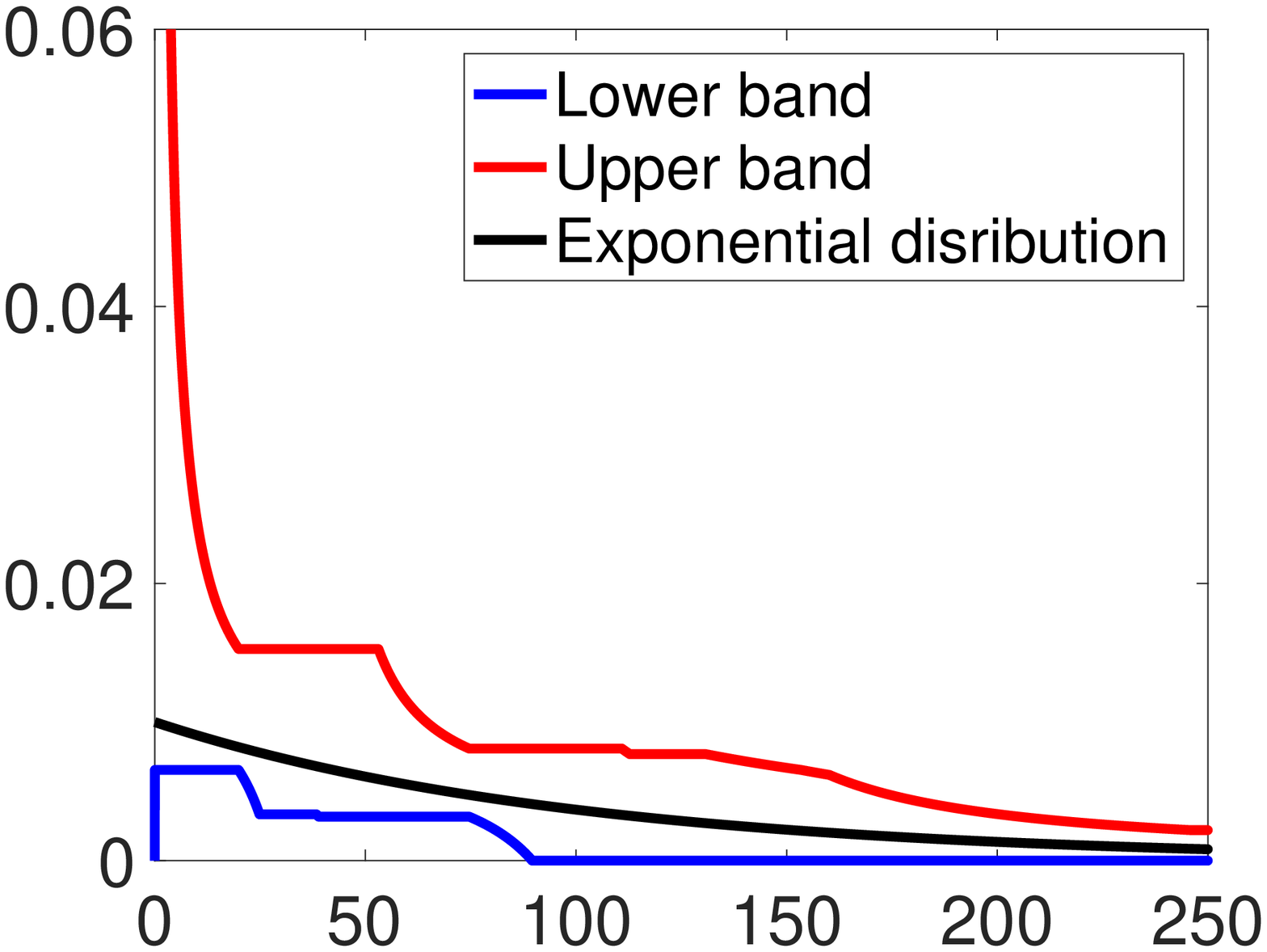}
			\caption{$\alpha=0.1$}
			\label{Fig:m}
		\end{subfigure}%
		\begin{subfigure}[b]{0.33\textwidth}
			\includegraphics[width=\linewidth]{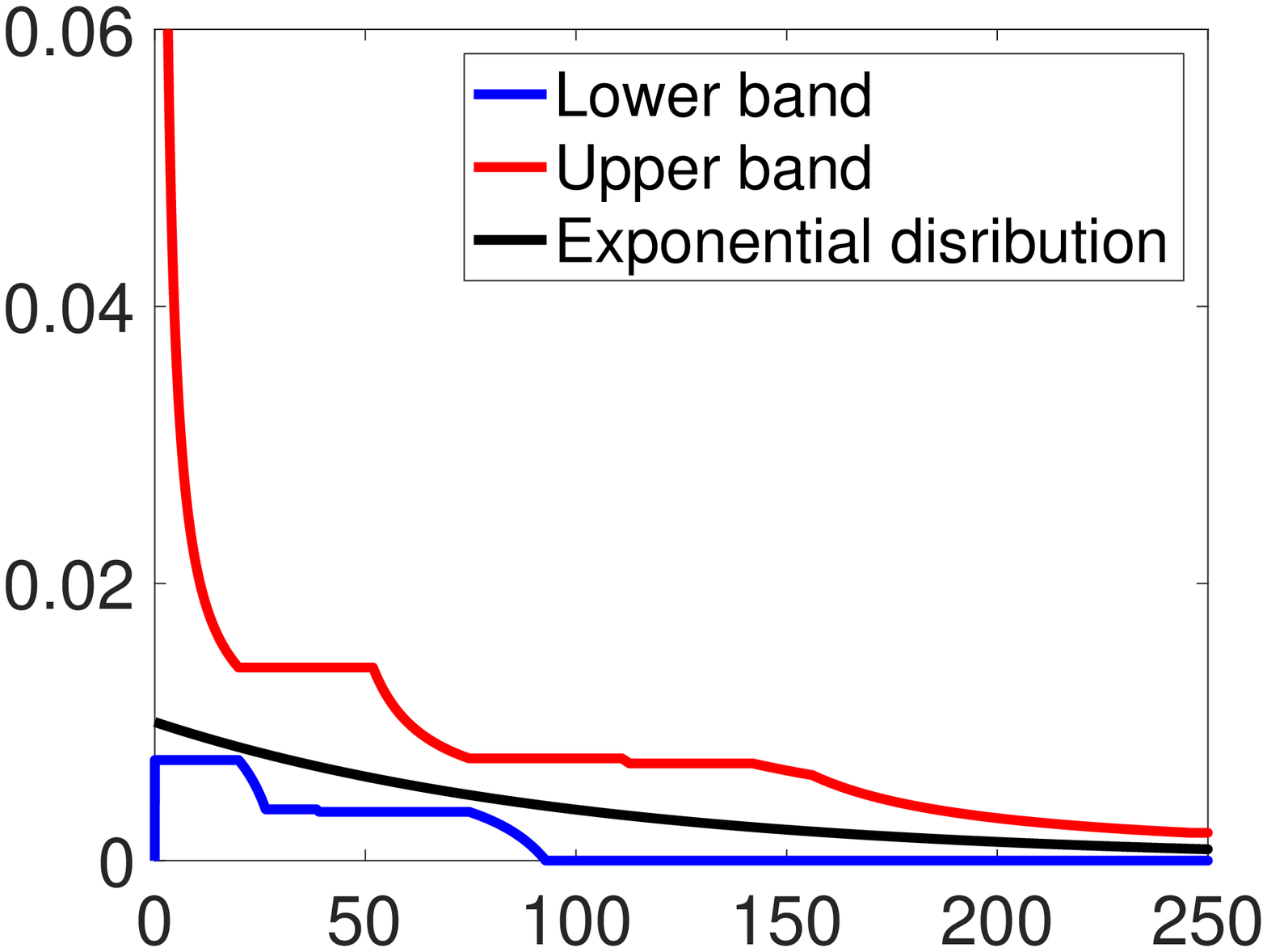}
			\caption{$\alpha=0.2$}
			\label{Fig:2m}
		\end{subfigure}%
		\begin{subfigure}[b]{0.33\textwidth}
			\includegraphics[width=\linewidth]{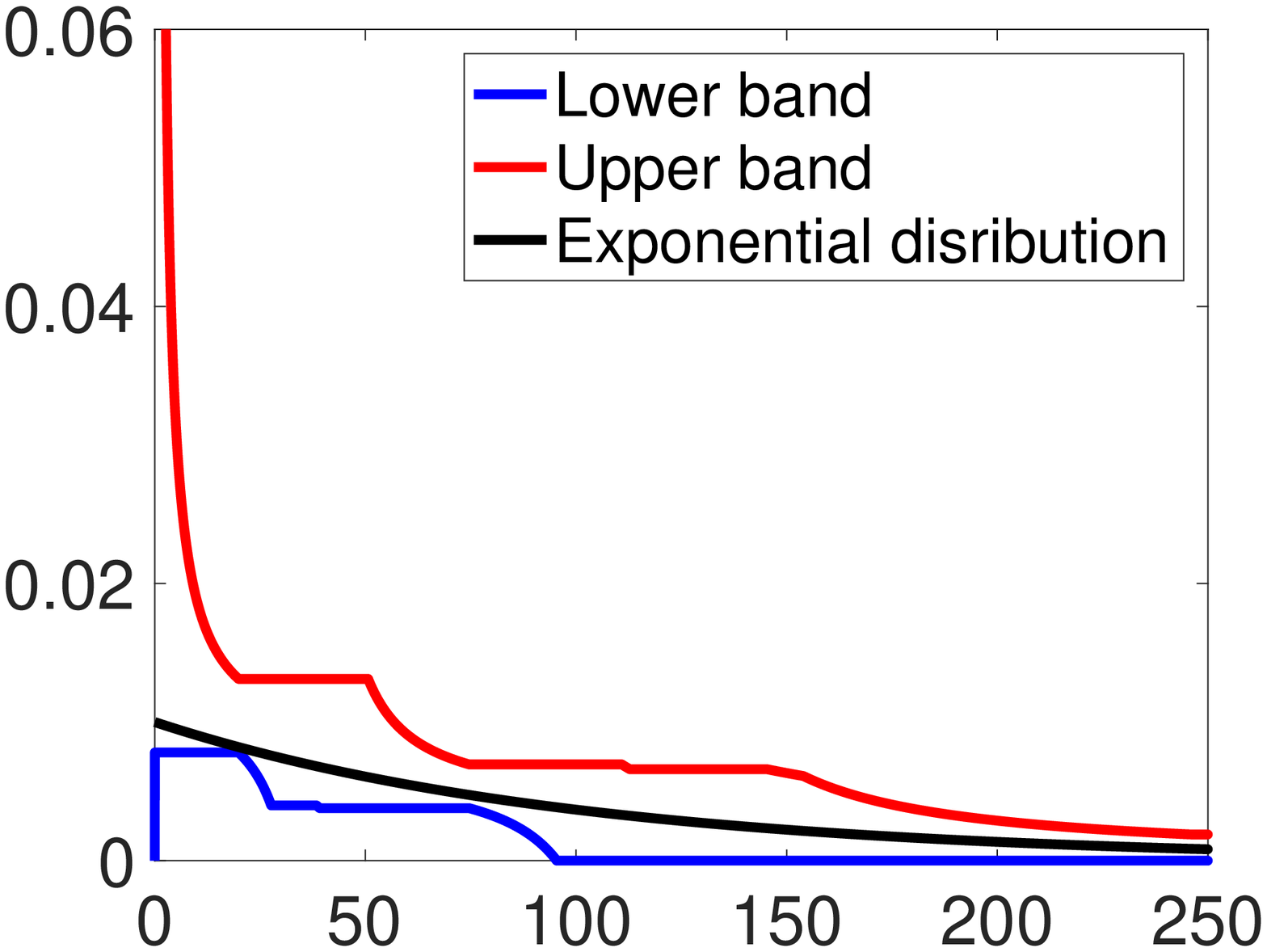}
			\caption{$\alpha=0.3$}
			\label{Fig:5m}
		\end{subfigure}%
		\quad
		\begin{subfigure}[b]{0.33\textwidth}
			\includegraphics[width=\linewidth]{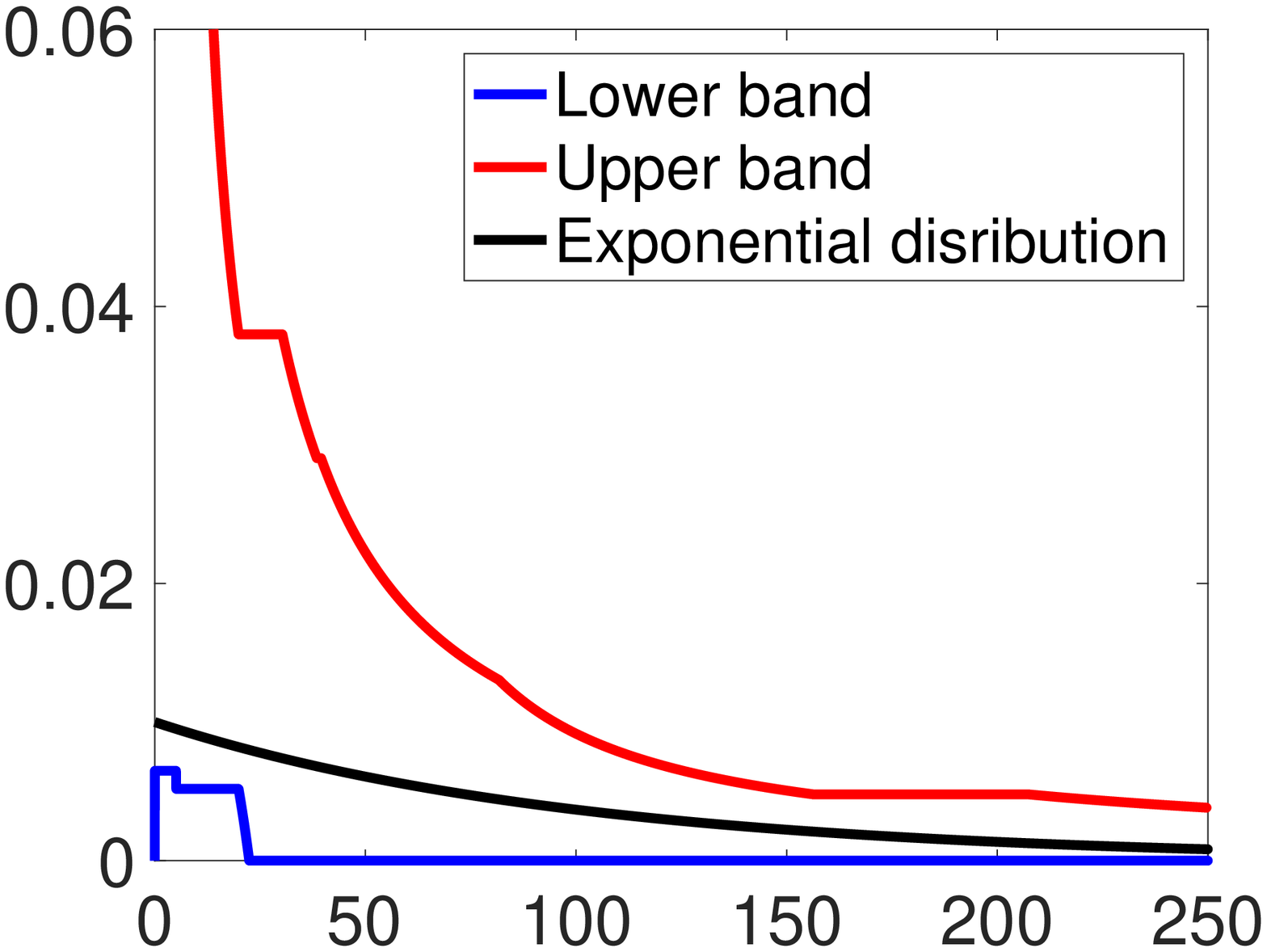}
			\caption{$N=10$}
			\label{Fig:m}
		\end{subfigure}%
		\begin{subfigure}[b]{0.33\textwidth}
			\includegraphics[width=\linewidth]{expo_100_02.eps}
			\caption{$N=100$}
			\label{Fig:2m}
		\end{subfigure}%
		\begin{subfigure}[b]{0.33\textwidth}
			\includegraphics[width=\linewidth]{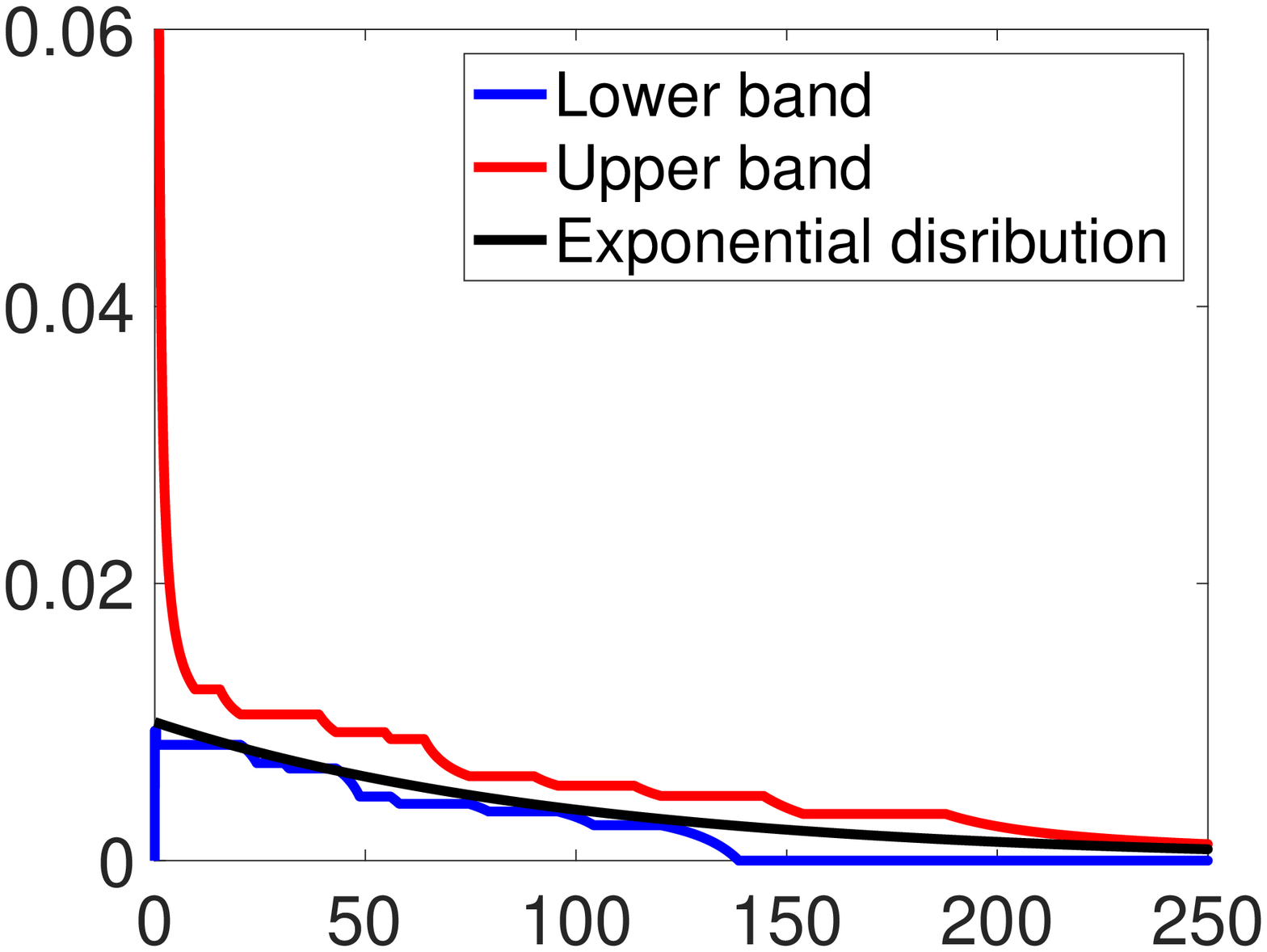}
			\caption{$N=1000$}
			\label{Fig:5m}
		\end{subfigure}
		\label{Fig:exp_mode}
	\end{figure}

	Finally, illustrations of this method are given in Figure \ref{Fig:beta_mode} 
	for a beta distribution and in Figure \ref{Fig:exp_mode} for a truncated 
	exponential distribution. (What exact parameters?)
	
	\subsection{Kernel-desnsity-estimation confidence bands}
	The shape-restricted confidence bands described in the previous section can only be applied to univariate density function. In this subsection, we describea method to construct confidence bands for multivariate density function based on the classical  \emph{kernel density estimation} (KDE)~\cite{rosenblatt1956remarks,parzen1962estimation}. This method requires a \emph{kernel function} which is a mapping $\mathcal{K}: \RR^m \rightarrow [0,+\infty)$ 
	satisfying $\int_{\RR^m} \mathcal{K}(\xi) d\xi = 1$. The commonly use kernel functions include uniform kernel
	$
	\mathcal{K}(\xi)=\frac{1}{2^m}\I_{\|\xi\|_{\infty}\leq 1}(\xi)
	$
	and Guassian kernel
	$
	\mathcal{K}(\xi)=\frac{1}{(2\pi)^{m/2}}\exp(-\|\xi\|_2^2/2)
	$.
	Let $h>0$ be a \emph{bandwidth parameter}. Recall that $\widehat \Xi_N := \{  \hat \xi_1, \ldots, \hat \xi_N\} \subseteq \RR^m$ is $N$ i.i.d. samples
	drawn from $p^*$. The KDE of $p^*$ based on $\widehat \Xi_N$, $\mathcal{K}$ and $h$ is 
	\begin{eqnarray}
	\label{eq:kerneldenstiy}
	\widehat p_h(\xi) := \dfrac{1}{N} \sum\limits_{i=1}^N \dfrac{1}{h^m}\mathcal{K}\left(\frac{\xi - \hat \xi_i}{h}\right).
	\end{eqnarray}
	The convergence of $\widehat p_h(\xi)$ to the true density $p^*(\xi)$ have been studied for a
	long time (see e.g. \cite{tsybakov2009introduction}) with most of existing works focusing the asymptotic convergence property. Recently, the finite-sample non-asymptotic convergence property of KDE is characterized by \cite{rinaldo2010generalized} and \cite{jiang2017uniform}. The confidence band we construct based 
	on KDE utilize the non-asymptotic convergence bound by \cite{jiang2017uniform}. 
	
	We need the following assumptions in this subsection.
	\begin{assumption}[For KDE confidence bands.]
		\label{assume:kernel}We assume:
		\begin{itemize}
			\item[A1.] There exists a constant $U$ such that $p(\mu) \leq U$ for any $\xi\in\Xi$.  
			\item[A2.] There exists a non-increasing function $\kappa:[0,+\infty)\rightarrow [0,+\infty)$ such that $\mathcal{K}(\xi)=\kappa(\|\xi\|_2)$.
			\item[A3.] There exists $r$, $C_r>0$ and $\tau>0$ such that for $t >\tau$, $\kappa(t)\leq C_r\cdot\exp(-t^r)$.
		\end{itemize}
	\end{assumption}
	A2 and A3 of Assumption~\ref{assume:kernel} hold if $\kappa$ is one of the popular kernel densities (up to scaling) in $\RR$, including the two mentioned above as well as Exponential, Tricube, triangular, and Epanechnikov kernels. Under these assumptions, the following finite-sample convergence result is established by \cite{jiang2017uniform}.

	\begin{proposition}\label{prop:KDEprop}[Theorem 2.~\cite{jiang2017uniform}] \label{lem:ratekde} Suppose $p^*$ is $(C,\rho)$-Holder continuous for constants $C>0$ and $\rho\in(0,1]$, i.e., 
		$\left|p^*(\xi')- p^*(\xi)\right| \leq C \left\|\xi' -\xi \right\|_2^{\rho}$ for any $\xi$ and $\xi'$.
		Let $\alpha\in(0,1)$ and $V_m$ be the volume of the unit ball in $\RR^m$. Suppose $h > \left(\frac{\log (N/\alpha)}{N}\right)^{1/m}$.
		We have 
		\small
		\begin{equation*}
		\label{eq:ratebandkde}
		\mathbb{P}\left\{\sup_{\xi\in\Xi}|\widehat p_h(\xi)-p^*(\xi)|\leq 
		C_1h^\rho+C_2\sqrt{\frac{\log(N/\alpha)}{N h^m}}
		\right\}\geq 1-\alpha.
		\end{equation*}
		\normalsize
		where\footnote{Note that $\int_0^{\infty}\kappa(t)t^{m+\rho}dt<+\infty$ because of A4 of Assumption~\ref{assume:kernel}.} $C_1=V_mC\int_0^{\infty}\kappa(t)t^{m+\rho}dt$ and $C_2=8m\sqrt{V_mU}\left(\int_0^{\infty}\kappa(t)t^{m/2}dt+1\right)+64m^2\kappa(0)$. As a consequence, if $h =\left(\frac{\log(N/\alpha)}{N}\right)^{1/(2\rho+m)}$, we have
		\small
		\begin{equation*}
		\label{eq:ratebandkdenew}
		\mathbb{P}\left\{\sup_{\xi\in\Xi}|\widehat p_h(\xi)-p^*(\xi)|\leq (C_1+C_2)\left(\frac{\log(N/\alpha)}{N}\right)^{\rho/(2\rho+m)}\right\}\geq 1-\alpha.
		\end{equation*}
		\normalsize
	\end{proposition}

	Based on the convergence property in Proposition~\ref{lem:ratekde}, with $h > \left(\frac{\log (N/\alpha)}{N}\right)^{1/m}$, we construct the KDE confidence band of a significance level of $\alpha$ for $p^*(\xi)$ as follows 
	\begin{equation} \label{equ:kernal_band}
	l_{\alpha}^{KDE} (\xi) = \max\{ 0, \ \widehat p_h(\xi) - \delta \}, \ u_{\alpha}^{KDE} (\xi) = \widehat p_h(\xi) + \delta
	\end{equation}
	where 
	\begin{equation} \label{equ:deltaKDE}
	\delta =  C_1h^\rho+C_2\sqrt{\frac{\log(N/\alpha)}{N h^m}}.
	\end{equation} 
	The corresponding uncertainty set is 
	\begin{eqnarray}
	\label{eq:PsetKDE}
	\DD^{KDE}(\widehat \Xi_N, \, \alpha) :=\left\{ p \in \LL \big| \begin{array}{c}
	l_{\alpha}^{KDE}(\xi) \le p(\xi) \le u_{\alpha}^{KDE}(\xi), \ \forall \, \xi \in  \Xi,  
	\int_{\Xi} p(\xi) \, d\xi = 1
	\end{array} \right\}.
	\end{eqnarray}
	The following property is a direct consequence of Proposition~\ref{lem:ratekde}. 
	\begin{theorem}
		\label{thm:KDEcover}
		$\PP \left\{p^*\in\mathcal{D}^{KDE}(\widehat \Xi_N,\alpha)\right\} \geq 1-\alpha$.
	\end{theorem}
	
	We summary this procedure for constructing a KDE confidence band in Algorithm~\ref{alg:KDE}. 
	\begin{algorithm}[h!]
		\caption{KDE-based confidence band $(l^{KDE}_\alpha(\xi), u^{KDE}_\alpha(\xi))$ at $\xi\in\Xi$}
		\label{alg:KDE}
		\begin{algorithmic}[1]
			\Require Data $\widehat \Xi_N = \{ \hat \xi^1, \ldots, \hat \xi^N \}$ sampled from $p^*$, a constant $U\geq p^*$ on $\Xi$, a significance level $\alpha\in(0,1)$, a kernel function $\mathcal{K}(\xi)=\kappa(\|\xi\|_2)$ with $\kappa$ satisfying Assumption~\ref{assume:kernel}, a bandwidth $h > \left(\frac{\log (N/\alpha)}{N}\right)^{1/m}$ and a point  $\xi\in\Xi$.	
			\State Compute $C_1=V_mC\int_0^{\infty}\kappa(t)t^{m+\rho}dt$ and $C_2=8m\sqrt{V_mU}\left(\int_0^{\infty}\kappa(t)t^{m/2}dt+1\right)+64m^2\kappa(0)$.
			\State Let $\widehat p_h(\xi)$ and $\delta$ defined as in \eqref{eq:kerneldenstiy} and \eqref{equ:deltaKDE}, respectively.
			\State Compute $(l_{\alpha}^{KDE}(\xi),u_{\alpha}^{KDE}(\xi))$ as in \eqref{equ:kernal_band}.
			\Ensure $(l_{\alpha}^{KDE}(\xi),u_{\alpha}^{KDE}(\xi))$
		\end{algorithmic}
	\end{algorithm}
	
	The following theorem show that, under additional assumption, $v_{\D^{KDE}(\widehat \Xi_N, \, \alpha)}^*$  define in \eqref{equ:aro} converges to $v^*$ in \eqref{eq:sp}.
	\begin{theorem} \label{thm:convergenceKDE}
		Suppose $p^*$ is $(C,\rho)$-Holder continuous for constants $C>0$ and $\rho>0$, i.e., 
		$\left|p^*(\xi')- p^*(\xi)\right| \leq C \left\|\xi' -\xi \right\|_2^{\rho}$ for any $\xi$ and $\xi'$. Moreover, suppose $\max_{x\in\mathcal{X},\xi\in \Xi}|f(x,\xi)|<+\infty$ and $h > \left(\frac{\log (N/\alpha)}{N}\right)^{1/m}$. For any $\epsilon>0$ and $\theta\in (0,1)$, there exists an $N_{\epsilon,\theta}$ such that for any $N\geq N_{\epsilon,\theta}$, we have
		\small
		\begin{eqnarray*}
			\PP\left(\sup_{x\in\mathcal{X}}\left|\int_{\Xi} f(x,\xi)p^*(\xi)d\xi-\sup_{p \in \DD^{KDE}(\widehat \Xi_N, \, \alpha)} \int_{\Xi} f(x,\xi) p(\xi) d\xi\right|\leq \epsilon\right)\geq 1-\theta
		\end{eqnarray*}
		\normalsize
		and
		$
		\PP\left(\left|v_{\D^{KDE}(\widehat \Xi_N, \, \alpha)}^*-v^*\right|\leq \epsilon\right)\geq 1-\theta.
		$
	\end{theorem}
	\proof{Proof.}
	See Appendix~\ref{sec:proof3}.
	%
	\endproof

	\section{An Numerical Method for DRO} \label{sec:method}
	
	In this section, we discuss a numerical scheme for solving the
	DRO problem in (\ref{equ:aro}) with an ambiguity set $\DD(\widehat \Xi_N,\alpha)$ in the form of \eqref{eq:Pset}. 
	Here, the ambiguity can be one of the two ambiguity sets $\DD^{SR}(\widehat \Xi_N,\alpha)$ and $\DD^{KDE}(\widehat \Xi_N,\alpha)$ introduced in Section~\ref{sec:ambiguityset} as long as Assumption~\ref{assume:Hengartner} or~\ref{assume:kernel} is satisfied. However, our numerical methods can be applied to a general ambiguity set like \eqref{eq:Pset} as long as the following assumptions are satisfied.
	We make the following assumption regarding to $\DD_{\alpha}(\widehat \Xi_N)$:
	\begin{assumption}[On ambiguity set $\DD(\widehat \Xi_N,\alpha)$] \label{ass:ul}~
		\begin{itemize}
			\item[a.] $\int |f(x,\xi)| u(\xi) d \xi < \infty, \quad \forall \, x \in \X$.
			\item[b.] $\DD(\widehat \Xi_N,\alpha) \neq \emptyset$.
		\end{itemize}
	\end{assumption}
	Note that, under Assumption~\ref{assume:Hengartner}, $\DD^{SR}(\widehat \Xi_N,\alpha)$ satisfies Assumption~\ref{ass:ul} because $u^{SR}_{\alpha}(\xi)\leq U$ for an upper bound $U$ of $p^*$ according to its definition~\eqref{eq:beta} and the constraints $0\leq\beta_j\leq U$ in \eqref{eq:lower_betac3}. Under Assumption~\ref{assume:kernel}, with the appropriate bandwidth $h$ (See Theorem~\ref{lem:ratekde}), we have $u_{\alpha}^{KDE} (\xi) = \widehat p_h(\xi) + \delta\leq p^*(\xi)+2\delta\leq U+2\delta$ with a high probability.
	
	
	We define $v(x)$ as the optimal value of the inner maximization
	problem of~\eqref{eq:Pset}, i.e.,
	\begin{equation} \label{equ:primal}
	v(x) := \sup_{p \; \in \; \DD(\widehat \Xi_N,\alpha)}   \int f(x,\xi) p(\xi) d\xi.
	\end{equation}
	By Assumption~\ref{ass:ul}, for every $x\in\X$,
	\begin{equation} \label{equ:primalfinite}
	|v(x)| \leq \sup_{p \; \in \; \DD(\widehat \Xi_N,\alpha)} \int |f(x,\xi)| p(\xi) d \xi \leq \int |f(x,\xi)| u(\xi) d \xi < \infty.
	\end{equation}
	
	Note that \eqref{equ:primal} is a \emph{continuous linear program} which are formulated with continuously many decision variable (or a functional decision variable) and continuously many constraints. In general, \eqref{equ:primal} cannot be reformulated as a convex optimization problem of finite dimension
	and solved by off-the-shelf optimization techniques as in most of the works in
	distributionally robust optimization. In this section, we propose a
	stochastic \emph{subgradient descent} (SGD) method for solving \eqref{equ:primal} and present its convergence rate.

	
	The dual problem of \eqref{equ:primal} is given as follows:
	\begin{equation}\label{equ:approx_dual3}
	\begin{array}{lll}
	 & \inf & \lambda -  \int l(\xi)\alpha(\xi)d\xi + \int u(\xi)\beta(\xi)d\xi   \\
	&\st & \lambda - \alpha(\xi) + \beta(\xi) \geq f(x, \xi) \ \ \ \forall \, \xi  \\
	& & \alpha(\xi) \geq 0, \ \beta(\xi) \geq 0 \ \ \ \forall \; \xi.
	\end{array}
	\end{equation}
	Weak duality always hold between (\ref{equ:primal}) and (\ref{equ:approx_dual3}).
	Furthermore, by the result of \cite{shapiro2001duality},
	strong duality holds between (\ref{equ:primal}) and (\ref{equ:approx_dual3}).
	Since $u(\xi)\geq l(\xi)\geq0$, it is easy to eliminate $\alpha(\xi)$ and $\beta(\xi)$
	in \eqref{equ:approx_dual3} and write it equivalently as:
	\begin{equation}\label{equ:approx_dual4}
	 v(x) =\inf  \limits_{\lambda} \lambda - \int l(\xi)(f(x, \xi) -\lambda)_-d\xi +\int u(\xi)(f(x, \xi) -\lambda)_+d\xi,
	\end{equation}
	where $(z)_-:=\max\{-z,0\}$ and $(z)_+:=\max\{z,0\}$. Thus, we have
	\begin{eqnarray} \label{equ:approx_dual5}
	v_{\DD(\widehat \Xi_N,\alpha)}^* :=\inf\limits_{x \in \X,\;\lambda}\left\{F(x,\lambda):=\lambda - \int l(\xi)(f(x, \xi) -\lambda)_-d\xi +\int u(\xi)(f(x, \xi) -\lambda)_+d\xi\right\}.
	\end{eqnarray}
In general, the two integrals appearing in \eqref{equ:approx_dual5} do not have an analytical form, raising a challenging for finding the optimal solution. In this section, we provide a stochastic gradient (SGD) method for solving \eqref{equ:approx_dual5} by approximating these integrals with random samples of $\xi$.
To apply gradient-based algorithm like SGD, we need to be able to compute the subgradient of $F(x,\lambda)$ in \eqref{equ:approx_dual5} with respect to both $x$ and $\lambda$. For doing that, we assume there exists a measurable mapping $f'(x,\xi):\mathcal{X}\times\Xi\rightarrow \mathbb{R}^d$ such that  $f'(x,\xi)\in\partial f(x,\xi)$ for any $x\in \mathcal{X}$ and $\xi\in\Xi$, where $\partial f(x,\xi)$ is the subdifferential of $f$ with respect to $x$. Then, under mild regularity conditions, the subgradients of   $F(x,\lambda)$ with respect to $x$ and $\lambda$ are, respectively, 
\begin{eqnarray}
\label{eq:deterministicg1}
\int l(\xi)f'(x,\xi)\I_{f(x, \xi) <\lambda}(\xi)d\xi+\int u(\xi)f'(x,\xi) \I_{f(x, \xi) \geq\lambda}(\xi)d\xi&\in&\partial_xF(x,\lambda)\\
\label{eq:deterministicg2}
1-\int l(\xi)\I_{f(x, \xi) <\lambda}(\xi)d\xi-\int u(\xi) \I_{f(x, \xi) \geq\lambda}(\xi)d\xi&\in&\partial_\lambda F(x,\lambda)
\end{eqnarray}	
where $\I$ is an indicator function. Then, we can use Monte Carlo method to approximate the integral in the subgradients above. 

In particular, let $I\subset\mathbb{R}^d$ be a set that contains the support of $l(\xi)$ and $u(\xi)$. Note that such a set always exists because $\Xi$ is assumed to be compact. We denote the volume of the box $I$ as $|I|$ and assume it can be computed. In fact, we can choose  $I\subset\mathbb{R}^d$ to be a box so that $|I|$ is the product of the lengths of all edges. Suppose $\xi$ is sampled from a uniform distribution on $I$. We can show that  \eqref{eq:deterministicg1} is the expectation of $|I|\I_{f(x, \xi) <\lambda}(\xi)l(\xi)f'(x, \xi)+|I|\I_{f(x, \xi) \geq\lambda}(\xi)u(\xi)f'(x, \xi)$ and  \eqref{eq:deterministicg2} is the expectation of $1-|I|\I_{f(x, \xi) <\lambda}(\xi)l(\xi)+|I|\I_{f(x, \xi) \geq\lambda}(\xi)u(\xi)$. Hence, we can use these two as the stochastic subgradient of $F$. Although the algorithm converges with any number of samples, we can apply the mini-batch techniques by generating $B$ i.i.d. samples from the uniform distribution on $I$ and constructing such a stochastic subgradient using each sample and then take the average of all samples. This mini-batch approach reduces the approximation noise of the Monte Carlo method and accelerate the algorithm in practice. 

Based on this idea, we proposed the SGD method for \eqref{equ:approx_dual5} in Algorithm~\ref{alg:SGD}. Note that $g_x$ and $g_\lambda$ are the mini-batch stochastic subgradients of $F$ with respect to $x$ and $\lambda$, respectively, with a batch size of $B$. The convergence analysis of Algorithm~\ref{alg:SGD} is standard and well-known (see e.g. \cite{nemirovski2009robust}) so we present the theorem below but omit its proof.

	\begin{algorithm}[h]
	\caption{SGD for \eqref{equ:approx_dual5}}\label{alg:SGD}
	\begin{algorithmic}[1]
		\Require An initial solution $(x_0,\lambda_0)\in\mathcal{X}\times\mathbb{R}$, batch size $B\geq 1$, step length $\eta_k=\eta/\sqrt{k+1}$ with $\eta>0$, a set $I\subset\mathbb{R}^d$ containing the support of $l$ and $u$, and the volume of $I$, denoted by $|I|$.
		\For{$k=0,1,\dots,$}
		\State $(\bar x_{k},\bar \lambda_{k})=\frac{\sum_{i=0}^{k}\eta_i(x_i,\lambda_i)}{\sum_{i=0}^{k}\eta_i}$
		\State Sample $\{\xi_1,\xi_2,\dots,\xi_B\}$ from a uniform distribution over $I$.
		\State Construct the stochastic gradients
		\begin{eqnarray*}
			g_x&=&\frac{|I|}{B}\sum_{i:f(x_k, \xi_i) <\lambda_k}l(\xi_i)f'(x_k, \xi_i)+\frac{|I|}{B}\sum_{i:f(x_k, \xi_i) \geq\lambda_k}u(\xi_i)f'(x_k, \xi_i)\\
			g_{\lambda}&=&1-\frac{|I|}{B}\sum_{i:f(x_k, \xi_i) <\lambda_k}l(\xi_i)-\frac{|I|}{B}\sum_{i:f(x_k, \xi_i) \geq\lambda_k}u(\xi_i)
		\end{eqnarray*}
		\State $x_{k+1}=\argmin_{x\in\mathcal{X}}\frac{1}{2}\|x-x_k+\eta_k g_x\|_2^2$ and $\lambda_{k+1}=\lambda_k-\eta_k g_{\lambda}$
		\EndFor
		\Ensure $(\bar x_{k},\bar \lambda_{k})$
	\end{algorithmic}
\end{algorithm}

\begin{theorem}
	\label{thm:converge}
Suppose there exists a constant $M$ such that $\mathbb{E}\|(g_x,g_\lambda)\|_2^2\leq M^2$. Algorithm~\ref{alg:SGD} ensures
\small
\begin{eqnarray}
\nonumber
\mathbb{E}\left[v(\bar x_k)-v_{\DD(\widehat \Xi_N,\alpha)}^*\right] \leq
\mathbb{E}\left[F(\bar x_k,\bar\lambda_k)-v_{\DD(\widehat \Xi_N,\alpha)}^*\right] &\leq& \frac{\|(x_0,\lambda_0)-(x_*,\lambda_*)\|_2^2+\sum_{i=0}^k\gamma_i^2M^2}{\sum_{i=0}^k\gamma_i}\\\label{eq:SGDconverge}
&\leq& \frac{\|(x_0,\lambda_0)-(x_*,\lambda_*)\|_2^2+\eta^2M^2(1+\ln(k+1))}{2\eta(\sqrt{k+2}-1))}
\end{eqnarray}
\normalsize
\end{theorem}
Note that the first inequality in \eqref{eq:SGDconverge} is because of \eqref{equ:approx_dual4} which indicates $v(x)\leq F(x,\lambda)$ for any $\lambda$.

\section{Computational Results} \label{sec:experiment}
	
	In this section, we validate our approach on two examples: a single-item
	newsvendor problem and a portfolio selection problem. Particularly,
	for the newsvendor example, we compare our approach with that in 
	\cite{bertsimas2014robust}, which applies hypothesis tests to construct
	ambiguity sets, and for the portfolio selection example, we compare our
	approach with that in \cite{esfahani2015data}, which applies the Wasserstein metric to construct ambiguity sets. 
	We implement Algorithm \ref{alg:SGD} in MATLAB 
	(R2014a) version 8.3.0.532. The linear programs from the approaches
	in the literature are solved by CPLEX 12.4 on an Intel Core i3 2.93 GHz
	Windows computer with 4GB of RAM and the mathematical models
	are implemented by using the modeling language YALMIP \cite{lofberg2004yalmip} 
	in MATLAB (R2014a) version 8.3.0.532.

	\subsection{Single-item newsvendor}
	
	We consider a classic single-item newsvendor problem in which we assume the demand
	$\xi \in \RR$ of an item follows a continuous distribution with a bounded support set $[a, b] \subseteq \RR$ with $0 \leq a < b$ and a bounded density function.
	An order of $x \geq 0$ units must be placed before demand occurs. After the 
	demand occurs, each unit of unmet demand incurs a shortage cost denoted by $c_s > 0$ 
	and each unit of surplus inventory incurs a holding cost denoted by $c_h > 0$. Hence, the cost function
	is defined as $f(x, \xi) = \max \left\{ c_s(\xi - x), \, c_h(x - \xi)\right\}$, which represents
	the cost of mismatch between supply and demand. In a classical newsvendor problem, 
	the goal is to determine the order size $x$ to minimize the expected cost. 
	When the demand's distribution is unknown, a corresponding DRO approach can be considered. Assuming a set of historical demand data is available,
	we construct an ambiguity set in the form of (\ref{eq:Pset}), or more specifically, $\DD^{SR}(\widehat \Xi_N, \, \alpha)$ in \eqref{eq:PsetSR} and solve the DRO (\ref{equ:aro}). We compare the optimal order obtained with the one found by the DRO model in \cite{bertsimas2014robust} where the ambiguity set is built using Kolmogorov-Smirnov test. 
	
	In our numerical experiments, we choose $c_s = 19$ and $c_h = 1$ and consider three different ground true distributions for the demand:  
	\begin{itemize}
		\item[1.] A truncated normal distribution created by truncating a normal distribution with mean 
		$100$ and standard deviation $50$ on $[0,250]$. 
		\item[2.] A beta distribution rescaled onto $[0,250]$ with parameters $\alpha=5$ and $\beta=2$. 
		\item[3.] A truncated exponential distribution created by truncating an exponential distribution with 
		mean $100$ on $[0,250]$.
	\end{itemize}
	
	For each distribution, we consider eight different sample sizes, i.e.,  
	$N \in \{10, 20, 40, 80\}$. For each size, we randomly generate a dataset $\widehat \Xi_N$ by i.i.d. sampling from the demand distribution. 	Using $\widehat \Xi_N$, we apply our approach and the method by \cite{bertsimas2014robust} to construct the ambiguity sets and then solve the corresponding  
	DRO problems to obtain an order size $\hat x$ from each approach. To evaluate the out-of-sample performance of $\hat x$, we sample
	another i.i.d. dataset $\{\xi'_i\}_{i=1}^{N_{\text{large}}}$ with $N_{\text{large}}=100,000$ from the true distribution and calculate the sample average approximation of the expected cost, i.e., $\frac{1}{N_{\text{large}}}\sum_{i=1}^{N_{\text{large}}}f(\hat x,\xi'_i)$ with $\hat x$ from each approach. We repeat this procedure 100 times to show the mean and variation of the  out-of-sample performance. 

When constructing the ambiguity set $\DD^{SR}(\widehat \Xi_N, \, \alpha)$ in our method, we need to provide a significance level $\alpha$ and a group size $K$ (see Algorithm~\ref{alg:cbHengartner}). Although a theoretical value of $K$ is suggested in Proposition~\ref{lem:rate}, it involves quantities which are hard to estimate (e.g. $B$ and $\rho$). When we construct $\DD^{SR}(\widehat \Xi_N, \, \alpha)$, we set $K=K(N,c) =\min\left\{\left\lceil c\left(N^{2}\log N\right)^{1/3}\right\rceil,N-1\right\}$ and select $c$ from $\{0.5,0.75,1,1.25,1.5\}$ based on the holdout validation method. The value significance level $\alpha$ is chosen from $\{0.75,0.8,0.85,0.95\}$  based on the same validation method. In particular, we randomly partition $\widehat \Xi_N$ into $\widehat \Xi_{train}$ and $\widehat \Xi_{text}$ with $|\widehat \Xi_{train}|=0.7N$ and $|\widehat \Xi_{test}|=0.3N$. Given a combination of $c$ and $\alpha$,  we first construct $\DD^{SR}(\widehat \Xi_{train}, \, \alpha)$ using and $K=K(0.7N,c)$, and then we solve $x^*_{c,\alpha} \in\argmin_{x \in \X} \sup_{p \in \DD^{SR}(\widehat \Xi_{train}, \, \alpha)} \int_{\Xi} f(x,\xi) p(\xi) d\xi$ using Algorithm~\ref{alg:SGD}. Then, we select the combination of $c$ and $\alpha$ with the largest  $\frac{1}{0.3N}\sum_{\hat\xi\in\widehat \Xi_{test}}f(x^*_{c,\alpha},\hat \xi)$. For each of the 100 independent trial, we repeat this process to select $c$ and $\alpha$. Note that the entire procedure is simulating how a decision maker select $K$ in practice when only $\widehat \Xi_N$ is available. For a fair comparison, we also apply the same validation scheme to choose the significance level $\alpha$ used in the method by \cite{bertsimas2014robust}.

\begin{figure} [t]
		\begin{center}
			\includegraphics[width=0.32\linewidth]{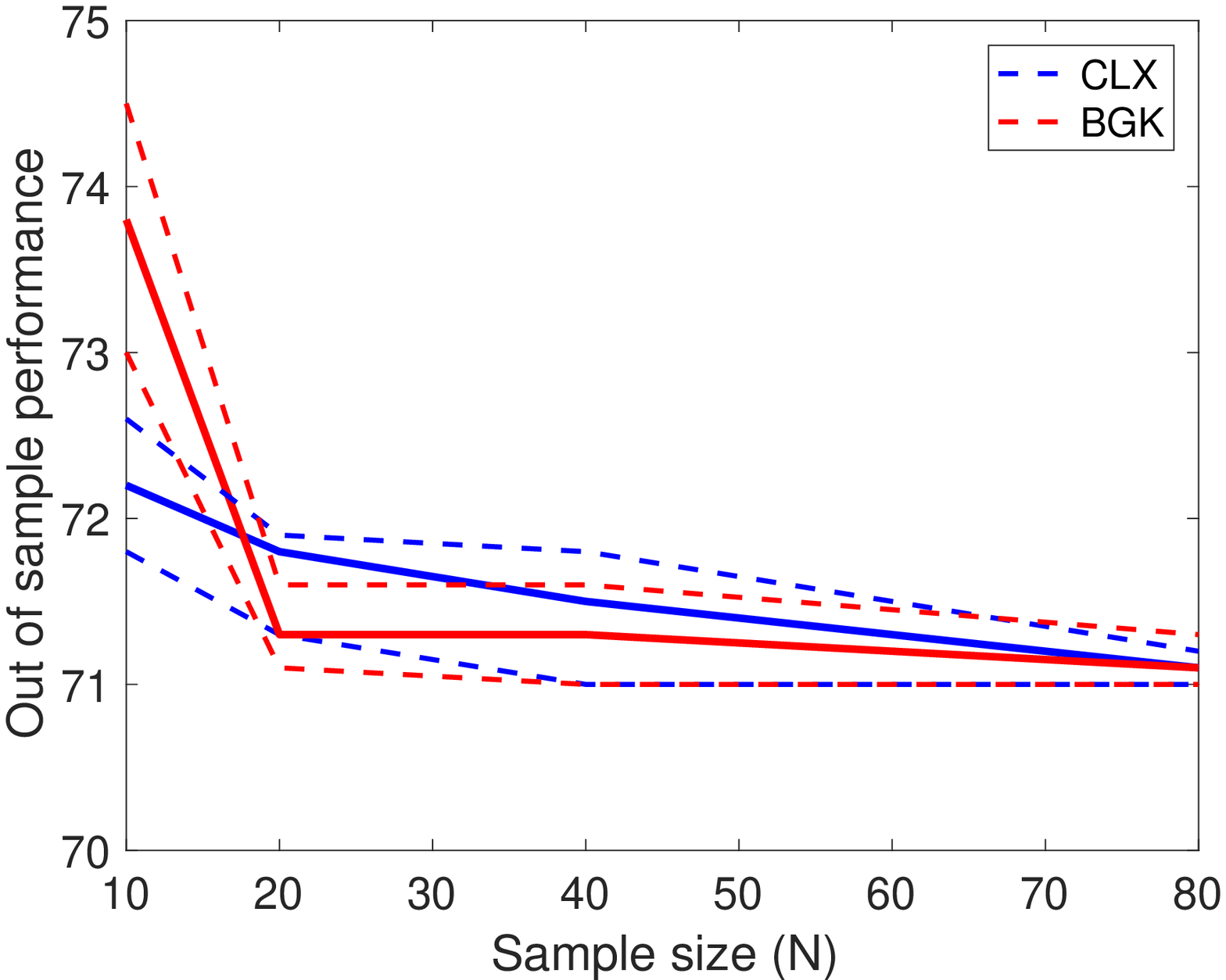}
			\includegraphics[width=0.32\linewidth]{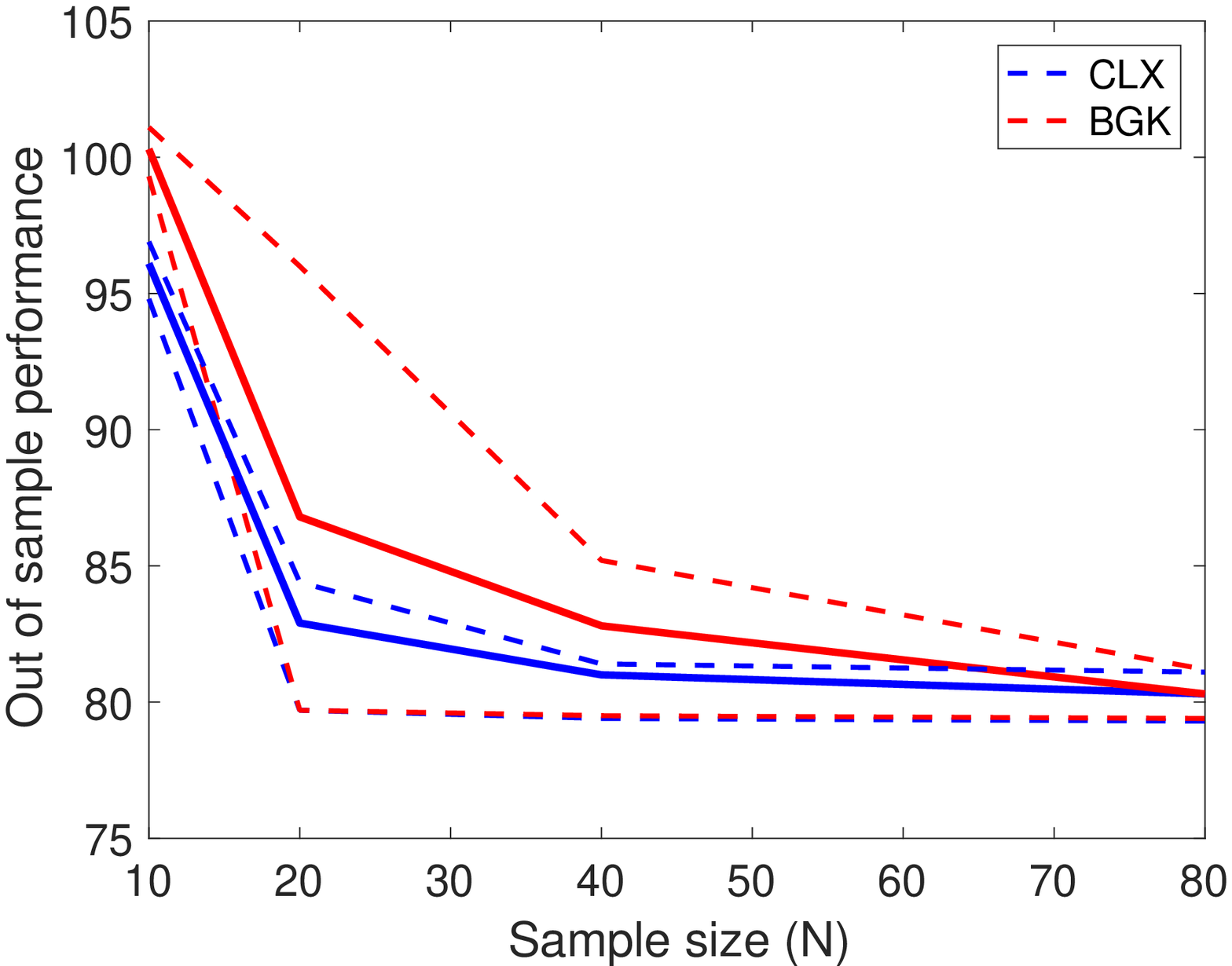}
			\includegraphics[width=0.32\linewidth]{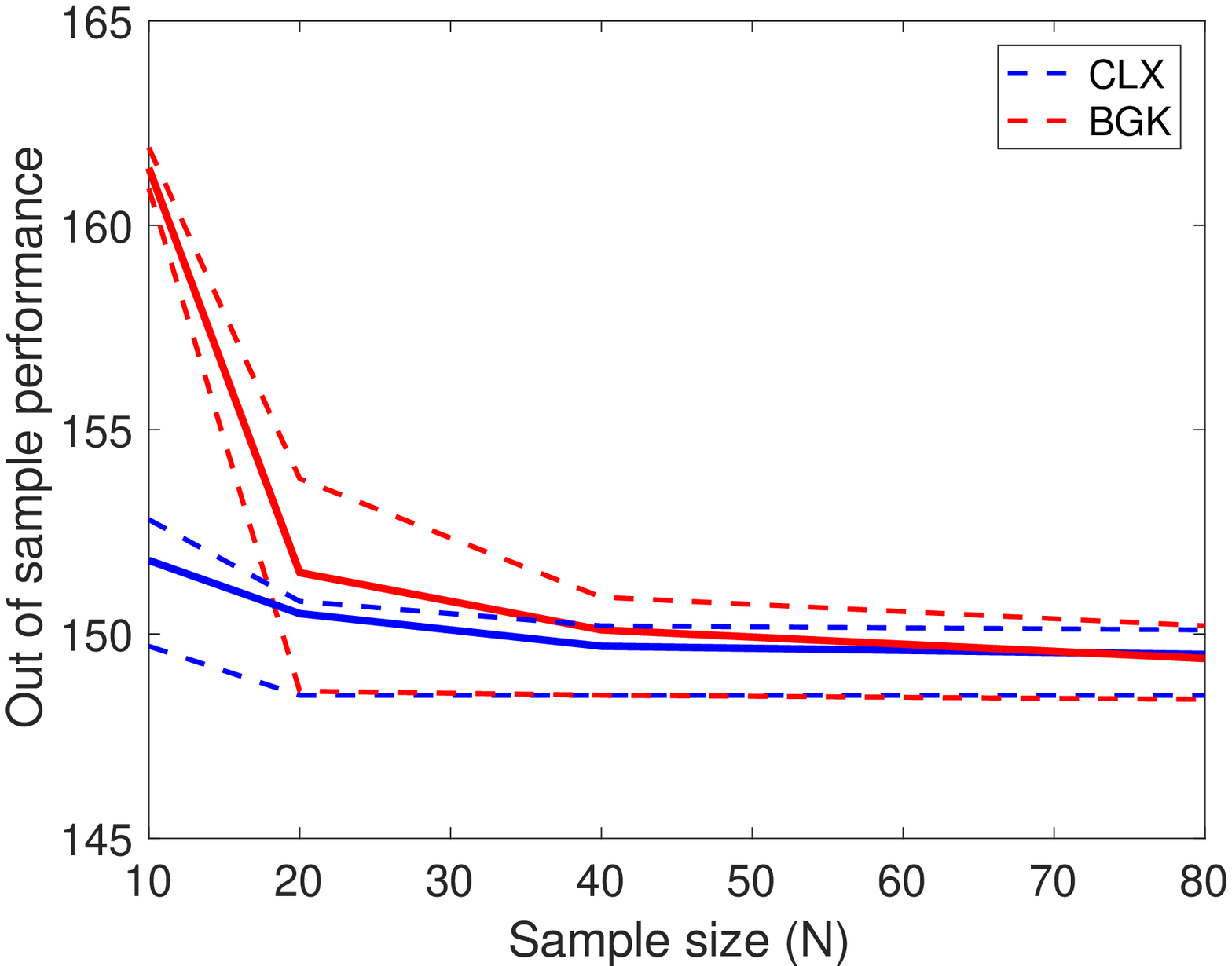}
				\caption{The out-of-sample performance of CLX and BGK on a newsvendor problem. The dash lines represent the $20^{\text{th}}$ and $80^{\text{th}}$ percentiles and the solid line represents the mean from 100 independent trials. Demand distribution: truncated normal (left), rescaled beta (middle), and truncated exponential (right).}
							\label{fig:single_norm}
		\end{center}
	\end{figure}

	We denote our approach by CLX and the method in \cite{bertsimas2014robust} by BGK. Figure \ref{fig:single_norm} illustrates the performances of CLX and BGK
	for each of the three distributions of demand. For each sample size $N$, we plot the $20^{\text{th}}$ percentile, the mean, and $80^{\text{th}}$ percentile of the out-of-sample performances in the 100 trials. The blue lines show the results from CLX, while the 
	red lines show the results from BGK. Figure \ref{fig:single_norm} indicates that CLX has a better out-of-sample performances than BGK when the sample size is small. As the sample size increases, the performances of both approaches become similar as the ambiguity sets in both approaches converge to the true deman distribution.

%
	
	
	
	\subsection{Portofolio management}
	
	In this example, we consider the classical portfolio selection problem consisting of  
	$n$ assets in which the investor must divide the total budget to 
	fractions $w=(w_1,w_2,\dots,w_n)$ 
	with $w_i\geq0$ and $\sum_{i=1}^nw_i=1$, and invest $w_i$ of the budget in 
	the $i^{\text{th}}$ security. We assume that the $i^{\text{th}}$ 
	security has a random future return $\xi_i$. The return from each unit of budget is 
	thus $w^\top\xi$. We assume that the investor is risk-averse and measures the investment risk by the \emph{conditional value at risk} (CVaR) of the return $w^\top\xi$; see \cite{rockafellar2000optimization}. Suppose the joint distribution of the return $\xi=(\xi_1,\xi_2,\dots,\xi_n)$ has a density function $p^*(\xi)$.
	The CVaR at level $\epsilon \in (0,1)$ of the return of a portfolio with respect to a probability distribution $p^*$ is defined as
	\[
	\text{CVaR}_{p^*,\epsilon}(-w^\top \xi) \equiv \inf \limits_{\beta \in \RR} \mathbb{E}_{p^*} \left[\beta + \frac{1}{\epsilon} \left( -w^\top\xi - \beta\right)_+\right]=\inf \limits_{\beta \in \RR} \int \left[\beta + \frac{1}{\epsilon} \left( -w^\top\xi - \beta\right)_+\right]p^*(\xi)d\xi,
	\]
	which represents the average of the $\epsilon \times 100\%$ worst portfolio
	losses (negative return) under distribution $p^*$. When $p^*$ is known, we consider the case where the investor wants to minimize a weighted sum of the mean and the
	CVaR of the portfolio loss $-w^\top \xi$, which is formulated as the following stochastic optimization:
	\begin{eqnarray*}
	&&\inf\limits_{w_i \geq 0, \sum_{i=1}^nw_i = 1} \left\{ \Ep_{p^*} [-w^\top \xi] +  \gamma \, \text{CVaR}_{p^*,\epsilon}(-w^\top \xi) \right\} \\
	& =& \inf\limits_{\beta,w_i \geq 0, \sum_{i=1}^nw_i = 1}  \int \left[ \max \left\{ -w^\top\xi +\gamma \beta, \, -(1+ \gamma/\epsilon)w^\top\xi + \gamma(1 - 1/\epsilon) \right\} \right]p^*(\xi)d\xi\\
	& =& \inf\limits_{x\in\mathcal{X}}  \int f(x,\xi)p(\xi)d\xi.
	\end{eqnarray*}
where $\gamma > 0$ indicates the investor's risk-aversion level, $\mathcal{X}=\{x=(w,\beta)\in\mathbb{R}^{n+1}|w_i \geq 0, \sum_{i=1}^nw_i = 1\}$, and 
$
f(x,\xi)=\max \left\{ -w^\top\xi +\gamma \beta, \, -(1+ \gamma/\epsilon)w^\top\xi + \gamma(1 - 1/\epsilon) \right\}.
$
	
	If the joint distribution $p^*$ of $\xi$ is unknown but a collection of 
	historical data return is collected, the investor can construct
	a data-driven ambiguity set of $p^*$ and solve the DRO problem corresponding to the stochastic optimization problem above. We construct the ambiguity set $\DD^{KDE}(\widehat \Xi_N, \, \alpha)$ in \eqref{eq:PsetKDE} and solve the DRO (\ref{equ:aro}) to construct an portfolio. Then, we compare our solution with the one obtained by the DRO model in \cite{esfahani2015data} where the ambiguity set is constructed using Wasserstein metric. 
	
	Following the numerical experiments by \cite{esfahani2015data}, we consider $n=10$ assets with decomposable returns 
	$\xi_i = \phi + \zeta_i$ for $i=1,2,\dots,10$ where $\phi \sim \text{normal}(0, 2\%)$ is a systematic
	risk factor shared by all assets and $\zeta_i \sim \text{normal}(i \times 3\%, 
	i \times 2.5\%)$ is an unsystematic risk factor associated with specific 
	assets. By the construction, assets with higher indices promise higher 
	mean returns at a higher risk. We set $\epsilon = 20\%$ and $\gamma = 10$
	in our all experiments. We consider $6$ different sample sizes, i.e.,  
	$N \in \{30, 60, 120, 240, 480, 960\}$. For each size, we randomly generate a dataset $\widehat \Xi_N$ by i.i.d. sampling returns from the aforementioned distribution of $\xi$. Using $\widehat \Xi_N$, we apply our approach and the method by \cite{esfahani2015data} to construct the ambiguity sets and then solve the DRO to obtain an portfolio $\hat x$ from each approach. To evaluate the out-of-sample performance of $\hat x$, we sample
	another i.i.d. dataset $\{\xi'_i\}_{i=1}^{N_{\text{large}}}$ with $N_{\text{large}}=100,000$ from the true distribution and calculate the sample average approximation of the expected cost, i.e., $\frac{1}{N_{\text{large}}}\sum_{i=1}^{N_{\text{large}}}f(\hat x,\xi'_i)$ with $\hat x$ from each approach. We repeat this procedure 100 times to show the mean and variation of the  out-of-sample performance. 
	
	When constructing $\DD^{KDE}(\widehat \Xi_N, \, \alpha)$, we choose $\mathcal{K}(\xi)=\kappa(\|\xi\|_2)$ wiht $\kappa$ being the boxcar kernel, namely, $\kappa(z)=Q$ if $z\in[0,1]$ and $\kappa(z)=0$ otherwise. Here, $Q$ is a normalization constant that ensures $\int_{\RR^m} \mathcal{K}(\xi) d\xi = 1$. Similar to $\DD^{SR}(\widehat \Xi_N, \, \alpha)$, the construction of $\DD^{KDE}(\widehat \Xi_N, \, \alpha)$ in our method requires some quantities which are hard to estimate (e.g. $C$ and $\rho$). Therefore, we construct  $\DD^{KDE}(\widehat \Xi_N, \, \alpha)$ using $l_{\alpha}^{KDE}$ and $u_{\alpha}^{KDE}$ in the form \eqref{equ:kernal_band} with the parameters $\delta$ and $h$ selected by the holdout validation method rather than their theoretical values in \eqref{equ:deltaKDE} and Proposition~\ref{prop:KDEprop}. In particular, we set  $h =c\left(\frac{\log(N)}{N}\right)^{1/(2+m)}$ (see Proposition~\ref{prop:KDEprop}) and then select $c$ from $\{0.02,0.04,0.06,0.08,0.1\}$ and $\delta$ from $\{0.02,0.04,0.06,0.08,0.1\}$. We randomly partition $\widehat \Xi_N$ into $\widehat \Xi_{train}$ and $\widehat \Xi_{text}$ with $|\widehat \Xi_{train}|=0.7N$ and $|\widehat \Xi_{test}|=0.3N$. Given a combination of $c$ and $\delta$,  we first construct $\DD^{KDE}(\widehat \Xi_{train})$ using and $K=K(0.7N,c)$, and then solve $x^*_{c,\delta} \in\argmin_{x \in \X} \sup_{p \in \DD^{KDE}(\widehat \Xi_{train})} \int_{\Xi} f(x,\xi) p(\xi) d\xi$ using Algorithm~\ref{alg:SGD}. Then, we select the combination of $c$ and $\delta$ with the largest  $\frac{1}{0.3N}\sum_{\hat\xi\in\widehat \Xi_{test}}f(x^*_{c,\delta},\hat \xi)$. For each of the 100 independent trial, we repeat this process to select $c$ and $\alpha$. For a fair comparison, we also apply the same validation scheme to choose the radius of the Wasserstein ball used to construct the ambiguity set in \cite{esfahani2015data}.

\begin{figure}[t]
		\begin{center}
			\includegraphics[width=0.7\linewidth]{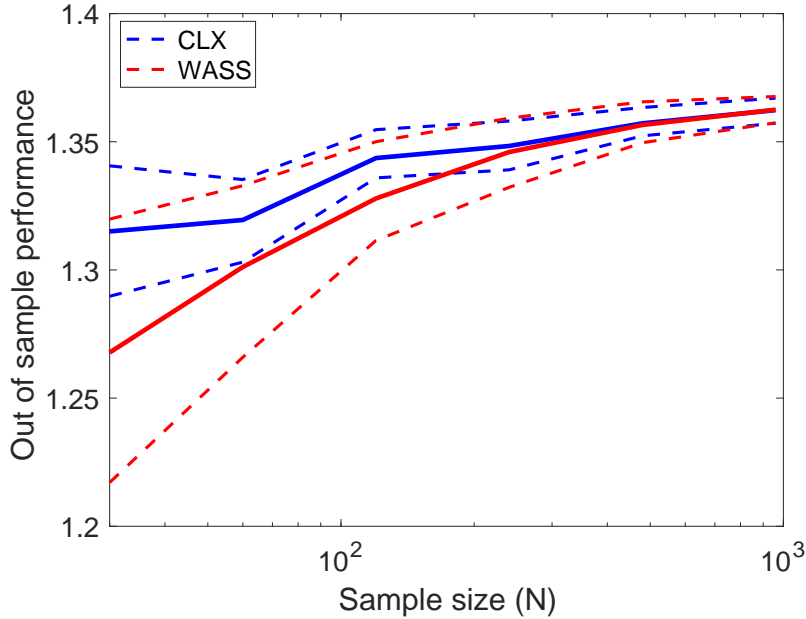}
			\caption{The out-of-sample performance of CLX and WASS on a portfolio selection problem. The dash lines represent the $20^{\text{th}}$ and $80^{\text{th}}$ percentiles and the solid line represents the mean from 100 independent trials. } \label{fig:portfolio}
		\end{center}
	\end{figure}

		We denote our approach by CLX and the method in \cite{esfahani2015data} by WASS and plot the numerical results in Figure \ref{fig:portfolio}.
	In particular, for each of the $100$ datasets in each sample-size scenario, 
	we evaluate the solutions from both CLX and Wasserstein by testing their 
	out-of-sample performances over the large dataset. We then plot the $20^{\text{th}}$ 
	percentile, the mean, and $80^{\text{th}}$ of the out-of-sample performances
	of both approaches over the $6$ sample-size scenarios. In figure 
	\ref{fig:portfolio}, the blue lines show the results from CLX, while the 
	red lines show the results from Wasserstein. Figure \ref{fig:portfolio} 
	indicates that both approaches converge to the true expectation with 
	the increase in sample size. 
	
	\section{Conclusions} \label{sec:conclude}
	
	In this paper, we proposed data-driven approaches to construct ambiguity 
	sets that consist of continuous probability density functions. The ambiguity 
	sets enjoy both finite and asymptotic convergences. The resulting distributionally
	robust optimization problem has infinite many variables and constraints. We
	then proposed a stochastic gradient decent method to solve the optimization
	problems. Numerical experiments in newsvendor problem and portfolio 
	selection problem verify the effectiveness of our approach.  
	
	%
	%
	%
	\section*{Appendix}
		
		
		\subsection*{Proof of Theorem~\ref{thm:convergence}}
		\label{sec:proof2}
		For $\delta>0$, we define a subset of $\Xi$ as
		$
		\Xi_\delta:=\{\xi\in\Xi|p^*(\xi)>\delta\}.
		$
		We first analyze the approximation error between $u_{\alpha}$, $l_{\alpha}$ and $p^*$ at a given point $\xi\in \Xi_\delta$. We assume $\xi>\mu$ first and the proof for $\xi<\mu$
		is similar. 
		For simplicity of notations, we use $c^{+}$ and $c^{-}$ to represent $c^{-}(\alpha)$ and $c^{-}(\alpha)$ in this proof. 
		
		There exists an index $i$ such that $\hat \xi_{(k_i)}< \xi< \hat \xi_{(k_{i+1})}$. Since $\xi>\mu$ and $p^*$ is $(C,\rho)$-Holder continuous, for a sufficiently large sample size $N$, we will further have $p^*( \hat \xi_{(k_{i+1})})\geq\frac{\delta}{2}$ (as $\hat \xi_{(k_{i+1})}$ is close enough to $\xi$) and $\mu<\hat \xi_{(k_{i-1})}<\hat \xi_{(k_i)}< \xi< \hat \xi_{(k_{i+1})}$. 
		Note that $p^*$ is monotonically decreasing over $[\hat \xi_{(k_i)}, \hat \xi_{(k_{i+1})}]$. By the definitions of $l_{\alpha}^{SR}(\xi)$ and $u_{\alpha}^{SR}(\xi)$ as in \eqref{eq:beta} and the constraints in \eqref{eq:lower_beta}, we must have $u\leq \frac{c^+}{\hat \xi_{(k_i)}-\hat \xi_{(k_{i-1})}}$ and $l\geq \frac{c^-}{\hat \xi_{(k_{i+1})}-\hat \xi_{(k_i)}}$, which implies
		\begin{eqnarray}
		\label{eq:betaandD}
		u_{\alpha}^{SR}(\xi)-l_{\alpha}^{SR}(\xi)\leq D_{\xi}\equiv\frac{c^+}{\hat \xi_{(k_i)}-\hat \xi_{(k_{i-1})}}-\frac{c^-}{\hat \xi_{(k_{i+1})}-\hat \xi_{(k_i)}}
		\end{eqnarray}
		
		Using the monotonicity of $p^*$, we can show that 
		\begin{eqnarray}
		\label{eq:introDelta1}
		(\hat \xi_{(k_i)}-\hat \xi_{(k_{i-1})})p^*(\hat \xi_{(k_i)})\leq\Delta_i &=& F^*(\hat \xi_{(k_i)})-F^*(\hat \xi_{(k_{(i-1)})})\leq (\hat \xi_{(k_i)}-\hat \xi_{(k_{i-1})})p^*(\hat \xi_{(k_{i-1})})\\
		\label{eq:introDelta2}
		(\hat \xi_{(k_{i+1})}-\hat \xi_{(k_i)})p^*(\hat \xi_{(k_{i+1})})\leq\Delta_{i+1} &=& F^*(\hat \xi_{(k_{i+1})})-F^*(\hat \xi_{(k_i)})\leq (\hat \xi_{(k_{i+1})}-\hat \xi_{(k_i)})p^*(\hat \xi_{(k_i)}).
		\end{eqnarray}
		Applying \eqref{eq:introDelta1} and \eqref{eq:introDelta2} to \eqref{eq:betaandD}, we obtain
		\small
		\begin{eqnarray}
		\nonumber
		D_{\xi}&\leq&\left|\frac{c^+p^*(\hat \xi_{(k_{i-1})})}{\Delta_i}-\frac{c^-p^*(\hat \xi_{(k_{i+1})})}{\Delta_{i+1}}\right|\\\nonumber
		&\leq&p^*(\hat \xi_{(k_{i+1})})\left|\frac{c^+}{\Delta_i}-\frac{c^-}{\Delta_{i+1}}\right|+\frac{c^+}{\Delta_i}\left|p^*(\hat \xi_{(k_{i+1})})-p^*(\hat \xi_{(k_{i-1})})\right|\\\nonumber
		&\leq&U\left|\frac{c^+}{\Delta_i}-\frac{c^-}{\Delta_{i+1}}\right|+\frac{Cc^+}{\Delta_i}\left|\hat \xi_{(k_{i+1})}-\hat \xi_{(k_{i-1})}\right|^\rho\\\nonumber
		&\leq&U\left|\frac{c^+}{\Delta_i}-\frac{c^-}{\Delta_{i+1}}\right|+\frac{Cc^+}{\Delta_i}\left|\frac{\Delta_i+\Delta_{i+1}}{p^*(\hat \xi_{(k_{i+1})})}\right|^\rho\\\label{eq:upperboundD}
		&\leq&U\left|\frac{c^+}{\Delta_i}-\frac{c^-}{\Delta_{i+1}}\right|+\frac{Cc^+}{\Delta_i}\left|\frac{\Delta_i+\Delta_{i+1}}{\delta/2}\right|^\rho
		\end{eqnarray}
		\normalsize
		
		According to equation (98) and (101) with ($\tau=2$) in ~\cite{Hengartner:95}, we have 
		\begin{eqnarray}
		\label{eq:probforD1}
		\lim_{N\rightarrow +\infty}\PP\left(\left|\frac{c^+}{\Delta_j}-\frac{c^-}{\Delta_{j+1}}\right|\leq 4\sqrt{\frac{\log(N/K)}{K}}\right)\geq 1-2\theta
		\end{eqnarray}
		and
		\begin{eqnarray}
		\label{eq:probforD2}
		\lim_{N\rightarrow +\infty}\PP\left(\frac{c^+}{\Delta_j}(\Delta_j+\Delta_{j+1})^\rho\leq 2\left(\frac{K}{N}\right)^\rho\right)\geq 1-2\theta
		\end{eqnarray}
		for any $\theta\in (0,1)$ and any $j=0,1,\dots,M$ (not necessarily $i$). 
		Recall that $\xi\in\Xi_\delta$ so that $p^*(\xi)>\delta$. We apply \eqref{eq:probforD1} and \eqref{eq:probforD2} to \eqref{eq:upperboundD} to achieve
		\begin{eqnarray}
		\label{eq:probforD3}
		\lim_{N\rightarrow +\infty}\PP\left(\max_{\xi\in\Xi_\delta}D_{\xi}\leq 4U\sqrt{\frac{\log(N/K)}{K}}+\frac{2^{\rho+1}C}{\delta^\rho}\left(\frac{K}{N}\right)^\rho\right)\geq 1-4\theta
		\end{eqnarray}
		
		Let $\mathcal{M}(\cdot)$ represent the Lebesgue measure on $\Xi$. We then choose $\delta$ to be small enough such that $2\max_{x\in\mathcal{X},\xi\in \Xi}|f(x,\xi)|U\mathcal{M}(\Xi\backslash\Xi_\delta)\leq\frac{\epsilon}{2}$.
		Consider a fixed $x\in\mathcal{X}$, we have 
		\begin{eqnarray}
		\nonumber
		&&\max_{x\in\mathcal{X}}\left|\int_{\Xi} f(x,\xi)p^*(\xi)d\xi-\int_{\Xi} f(x,\xi)u_\alpha^{SR}(\xi)d\xi\right|\\\nonumber
		&\leq&\max_{x\in\mathcal{X}}\left|\int_{\Xi_\delta} f(x,\xi)p^*(\xi)d\xi-\int_{\Xi_\delta} f(x,\xi)u_\alpha^{SR}(\xi)d\xi\right|+\max_{x\in\mathcal{X}}\left|\int_{\Xi\backslash\Xi_\delta} f(x,\xi)p^*(\xi)d\xi-\int_{\Xi\backslash\Xi_\delta} f(x,\xi)u_\alpha^{SR}(\xi)d\xi\right|\\\nonumber
		&\leq&\max_{x\in\mathcal{X},\xi\in \Xi}|f(x,\xi)|\max_{\xi\in\Xi_\delta}D_{\xi}+2\max_{x\in\mathcal{X},\xi\in \Xi}|f(x,\xi)|U\mathcal{M}(\Xi\backslash\Xi_\delta)
		\end{eqnarray}
		which, according to \eqref{eq:probforD3}, implies
		\small
		\begin{eqnarray*}
			\lim_{N\rightarrow +\infty}\PP\left(
			\begin{array}{ll}
				&\max_{x\in\mathcal{X}}\left|\int_{\Xi} f(x,\xi)p^*(\xi)d\xi-\int_{\Xi} f(x,\xi)u_\alpha(\xi)d\xi\right|\\
				\leq& \max_{x\in\mathcal{X},\xi\in \Xi}|f(x,\xi)|\left[4U\sqrt{\frac{\log(N/K)}{K}}+\frac{2^{\rho+1}C}{\delta^\rho}\left(\frac{K}{N}\right)^\rho\right]+\frac{\epsilon}{2}
			\end{array}
			\right)\geq 1-4\theta
		\end{eqnarray*}
		\normalsize
		Similarly, we can also show that 
		\small
		\begin{eqnarray*}
			\lim_{N\rightarrow +\infty}\PP\left(
			\begin{array}{ll}
				&\max_{x\in\mathcal{X}}\left|\int_{\Xi} f(x,\xi)p^*(\xi)d\xi-\int_{\Xi} f(x,\xi)l_\alpha(\xi)d\xi\right|\\
				\leq& \max_{x\in\mathcal{X},\xi\in \Xi}|f(x,\xi)|\left[4U\sqrt{\frac{\log(N/K)}{K}}+\frac{2^{\rho+1}C}{\delta^\rho}\left(\frac{K}{N}\right)^\rho\right]+\frac{\epsilon}{2}
			\end{array}
			\right)\geq 1-4\theta
		\end{eqnarray*}
		\normalsize
		Therefore, for any $\epsilon$ and $\theta$, there exists an $N_{\epsilon,\theta}$ such tha,t for any $N\geq N_{\epsilon,\theta}$, 
		\small
		\begin{eqnarray*}
			\PP\left(\max_{x\in\mathcal{X}}\left|\int_{\Xi} f(x,\xi)p^*(\xi)d\xi-\sup_{p \in \DD^{SR}(\widehat \Xi_N, \, \alpha)} \int_{\Xi} f(x,\xi) p(\xi) d\xi\right|\leq \epsilon\right)\geq 1-10\theta.
		\end{eqnarray*}
		\normalsize
		Then, we have proof the first conclusion. The second conclusion can be easily implied from the first conclusion.
		
		\subsection*{Proof of Theorem~\ref{thm:convergenceKDE}}
		\label{sec:proof3}
		Suppose $p^*\in\DD^{KDE}(\widehat \Xi_N, \, \alpha)$. We have
		\begin{eqnarray*}
			\nonumber
			&&\max_{x\in\mathcal{X}}\left|\int_{\Xi} f(x,\xi)p^*(\xi)d\xi-\sup_{p \in \DD^{KDE}(\widehat \Xi_N, \, \alpha)}\int_{\Xi} f(x,\xi)p(\xi)d\xi\right|\\
			&\leq&\max_{x\in\mathcal{X}}\int_{\Xi} |f(x,\xi)||u_{\alpha}^{KDE}(\xi)-l_{\alpha}^{KDE}(\xi)|d\xi\\
			&\leq& \max_{x\in\mathcal{X}}\int_{\Xi} |f(x,\xi)|2\delta d\xi\leq2\delta\max_{x\in\mathcal{X},\xi\in\Xi}|f(x,\xi)|\mathcal{M}(\Xi)
		\end{eqnarray*} 
		where $\delta$ is defined in \eqref{equ:deltaKDE}. By the definition of $\delta$, there exists $N_{\epsilon,\alpha}$ such that for $N>N_{\epsilon,\alpha}$, we have
		$2\delta\max_{x\in\mathcal{X},\xi\in\Xi}|f(x,\xi)|\mathcal{M}(\Xi)\leq \epsilon$. As a result of Theorem~\ref{thm:KDEcover}, for $N>N_{\epsilon,\alpha}$, we have
		\begin{eqnarray*}
			&&\PP\left(\max_{x\in\mathcal{X}}\left|\int_{\Xi} f(x,\xi)p^*(\xi)d\xi-\sup_{p \in \DD^{KDE}(\widehat \Xi_N, \, \alpha)}\int_{\Xi} f(x,\xi)p(\xi)d\xi\right|\leq \epsilon\right)\\
			&\geq& \PP\left(p^*\in\DD^{KDE}(\widehat \Xi_N, \, \alpha)\right)\geq 1-\alpha
		\end{eqnarray*} 
		The second conclusion can be easily implied from the first conclusion.

	\section*{Acknowledgements}
	The authors would like to thank Aditya
		Guntuboyina for referring us to the papers about confidence band constructions.

\end{onehalfspace}

\bibliographystyle{plain}

\bibliography{dro}

\appendix

\end{document}